\documentclass[11pt]{article}

\usepackage{graphicx}
\usepackage[colorlinks=true,breaklinks=true,bookmarks=true,urlcolor=blue,
     citecolor=blue,linkcolor=blue,bookmarksopen=false,draft=false]{hyperref}
\usepackage[mathscr]{eucal}
\usepackage{epsfig,epsf,psfrag}
\usepackage{amssymb,amsfonts,latexsym}
\usepackage{amsmath}
\usepackage{graphicx}
\usepackage{epstopdf}
\usepackage{bm,xcolor,url}
\usepackage{fixltx2e}%ordering of single and double column floats
\usepackage{array}%array and tabular environments
\usepackage{verbatim}
\usepackage{algpseudocode, cite}
\usepackage{algorithm}
\usepackage{verbatim}
\usepackage{textcomp}
\usepackage{mathrsfs}
\usepackage[referable]{threeparttablex}
\usepackage{amsthm}
\usepackage{enumerate}
\usepackage{footnote}
\usepackage{enumitem}
\usepackage{bbm}
\usepackage{xr}
\usepackage{mathtools}
\catcode`~=11 \def\UrlSpecials{\do\~{\kern -.15em\lower .7ex\hbox{~}\kern .04em}} \catcode`~=13 

\allowdisplaybreaks[3]

\newcommand{\tnorm}[1]{{\left\vert\kern-0.25ex\left\vert\kern-0.25ex\left\vert #1 
    \right\vert\kern-0.25ex\right\vert\kern-0.25ex\right\vert}}
\newcommand{\tnormt}[1]{{\vert\kern-0.25ex\vert\kern-0.25ex\vert #1 
    \vert\kern-0.25ex\vert\kern-0.25ex\vert}}

\newcommand{\normt}[1]{\Vert#1\Vert}
\newcommand{\abst}[1]{\vert#1\vert}
\newcommand{\abs}[1]{\left\lvert#1\right\rvert}

\newcommand{\nn}{\nonumber}

\newcommand{\nt}{\addtocounter{equation}{1}\tag{\theequation}} 

\newcommand{\dom}{\mathsf{dom}\,}

\newcommand{\inter}{\mathsf{int}\,}
\newcommand{\bdry}{\mathsf{bd}\,}

\newcommand{\ri}{\mathsf{ri}\,}
\newcommand{\aff}{\mathsf{aff}\,}
\newcommand{\lin}{\mathsf{lin}\,}

\newcommand{\conv}{\mathsf{conv}\,}
\newcommand{\cone}{\mathrm{cone}\,}

\newcommand{\Diam}{\mathsf{Diam}\,}
\renewcommand{\dim}{\mathsf{dim}\,}
\newcommand{\dist}{\mathsf{dist}\,}

\newcommand{\bbone}{\mathbbm{1}}

\newcommand{\where}{\;\; \mbox{where}\;\;}
\newcommand{\andd}{\;\; \mbox{and}\;\;}
\newcommand{\Ext}{\mathsf{Ext}}

\newcommand{\baralpha}{\bar{\alpha}}

\newcommand{\calA}{\mathcal{A}}
\newcommand{\calB}{\mathcal{B}}

\newcommand{\calD}{\mathcal{D}}

\newcommand{\calF}{\mathcal{F}}
\newcommand{\calG}{\mathcal{G}}
\newcommand{\calH}{\mathcal{H}}
\newcommand{\calI}{\mathcal{I}}
\newcommand{\calJ}{\mathcal{J}}
\newcommand{\calK}{\mathcal{K}}

\newcommand{\calN}{\mathcal{N}}

\newcommand{\calS}{\mathcal{S}}
\newcommand{\calT}{\mathcal{T}}

\newcommand{\calV}{\mathcal{V}}
\newcommand{\calW}{\mathcal{W}}
\newcommand{\calX}{\mathcal{X}}
\newcommand{\calY}{\mathcal{Y}}
\newcommand{\calZ}{\mathcal{Z}}

\newcommand{\barcalX}{\bar{\calX}}

\newcommand{\rmA}{\mathrm{A}}

\newcommand{\rmd}{\mathrm{d}}

\newcommand{\rmF}{\mathrm{F}}

\newcommand{\bbR}{\mathbb{R}}
\newcommand{\bbS}{\mathbb{S}}

\newcommand{\bbX}{\mathbb{X}}
\newcommand{\bbY}{\mathbb{Y}}
\newcommand{\bbZ}{\mathbb{Z}}

\DeclareMathAlphabet{\mathbsf}{OT1}{cmss}{bx}{n}
\newcommand{\rvA}{\mathsf{A}}

\newcommand{\rvE}{\mathsf{E}}

\newcommand{\tilG}{\widetilde{G}}

\newcommand{\barn}{\bar{n}}

\newcommand{\barp}{\bar{p}}

\newcommand{\barv}{\overline{v}}
\newcommand{\barw}{\overline{w}}

\newcommand{\bary}{\bar{y}}

\newcommand{\barF}{\bar{F}}

\newcommand{\ceil}[1]{\lceil{#1}\rceil}
\newcommand{\floor}[1]{\lfloor{#1}\rfloor}

\newcommand{\lranglet}[2]{\langle{#1},{#2}\rangle}

\newcommand{\ipt}[2]{\langle{#1},{#2}\rangle}

\newcommand{\lea}{\stackrel{\rm(a)}{\le}}

\newcommand{\gea}{\stackrel{\rm(a)}{\ge}}
\newcommand{\geb}{\stackrel{\rm(b)}{\ge}}

\DeclareMathOperator*{\argmax}{arg\,max}
\DeclareMathOperator*{\argmin}{arg\,min}

\DeclareMathOperator{\st}{s.t.}

\DeclareMathOperator{\supp}{\mathsf{supp}}

\newtheorem{theorem}{Theorem} 
\newtheorem*{theorem*}{Theorem}
\newtheorem{lemma}{Lemma}

\newtheorem*{assump*}{Assumption}

\theoremstyle{remark}
\newtheorem{remark}{Remark}

\newcommand{\qednew}{\nobreak \ifvmode \relax \else
      \ifdim\lastskip<1.5em \hskip-\lastskip
      \hskip1.5em plus0em minus0.5em \fi \nobreak
      \vrule height0.75em width0.5em depth0.25em\fi}

\usepackage{fullpage}
\numberwithin{equation}{section}
\numberwithin{lemma}{section}
\numberwithin{assump}{section}
\numberwithin{theorem}{section}
\numberwithin{prop}{section}
\numberwithin{remark}{section}
\numberwithin{corollary}{section}
\newcommand{\poi}{(P)}

\newcommand{\barDelta}{\bar\Delta}
\usepackage[labelformat=simple]{subcaption}

\newcommand{\normtt}[1]{{\left\vert\kern-0.25ex\left\vert\kern-0.25ex\left\vert\kern 0.25ex  #1 \kern 0.25ex    \right\vert\kern-0.25ex\right\vert\kern-0.25ex\right\vert}}

\begin{document}

\title{New Analysis of An Away-Step Frank-Wolfe Method for Minimizing Logarithmically-Homogeneous Barriers} %: Global and Local Linear Convergence % and Conditioning of D-Optimal Design %and Positron Emission Tomography   %, {with an Extension}
%\thanks{Research
%supported by AFOSR Grant No. FA9550-19-1-0240.



\author{Renbo Zhao\thanks{Department of Business Analytics, Tippie College of Business, University of Iowa (\href{mailto:renbo-zhao@uiowa.edu}{renbo-zhao@uiowa.edu}).}
      %  Robert M. Freund\thanks{MIT Sloan School of Management, 77 Massachusetts Avenue, Cambridge, MA   02139 ({mailto:  rfreund@mit.edu}).} %etc.
}


%\institute{Renbo Zhao \at
%			  Operations Research Center\\
%              Massachusetts Institute of Technology\\
%              \email{renboz@mit.edu}           %  \\
%%             \emph{Present address:} of F. Author  %  if needed
%           Robert M. Freund \at
%           	  Sloan School of Management\\
%              Massachusetts Institute of Technology\\
%              \email{rfreund@mit.edu}

\maketitle

\begin{abstract} %\textcolor{red}{Abstract is TBD.}

We present and analyze an away-step Frank-Wolfe method for the convex optimization problem ${\min}_{x\in\calX}  \; f(\mathsf{A} x) + \ipt{c}{x}$, where $f$ is a $\theta$-logarithmically-homogeneous self-concordant barrier, $\rvA$ is a linear operator that may be non-invertible, $\ipt{c}{\cdot}$ is a linear function and $\calX$ is a nonempty polytope. %This class of problems include three 
The applications of primary interest %that we are primarily interested in 
include D-optimal design, inference of multivariate Hawkes processes, and TV-regularized Poisson image de-blurring. %with total-variation regularization. 
We establish affine-invariant and norm-independent global linear convergence rates of our method, in terms of both the objective gap and the Frank-Wolfe gap. When specialized to the D-optimal design problem, our results settle a question left open since Ahipasaoglu, Sun and Todd~\cite{Ahi_08}. We also show that the iterates generated by our method will land on and remain in a face of $\calX$ within a bounded number of iterations, which can lead to improved local linear convergence rates (for both the objective gap and the Frank-Wolfe gap). We conduct 
numerical experiments on %the problems of 
D-optimal design and inference of multivariate Hawkes processes, and our  results  not only  demonstrate the efficiency and effectiveness of our method compared to other principled first-order methods, but also corroborate our theoretical results quite well. 

%raised in Ahipasaoglu, Sun and Todd (2008) on the global linear convergence of the away-step Frank-Wolfe method specialized to the D-optimal design problem. 

\end{abstract}

%\noindent {\bf Keywords:} \vspace{0.2in}

\section{Introduction}\label{sec:intro}

We present and analyze an away-step Frank-Wolfe (FW) method~\cite{Wolfe_70,Guelat_86,Lacoste_15,Beck_17,Pena_19,Garber_20} for the following class of optimization problems:
\begin{equation}
 F^*:= {\min}_{x\in\calX} \;[F(x):= f(\rvA x) + \lranglet{c}{x}]. \tag{P} \label{eq:poi}
\end{equation}
In~\eqref{eq:poi},  $\rvA:\bbX\to \bbY$ is a (potentially) non-invertible linear operator for some finite-dimensional vector spaces $\bbX$ and $\bbY$, and  $c\in\bbX^*$ is a linear function on $\bbX$, where $\bbX^*$ denotes the dual space of $\bbX$. In addition, $f:\inter\calK\to\bbR$ is a $\theta$-logarithmically-homogeneous self-concordant barrier ($\theta$-LHSCB) on some regular cone $\calK\subseteq \bbY$, where $\inter\calK$ denotes the interior of $\calK$. (The definition of $\theta$-LHSCB will be reviewed in Section~\ref{sec:LHSCB}.)  % that is differentiable on its domain $\dom f:=\{x \in \mathbb{R}^m: f(x) < +\infty \}$,  and the linear functional $c\in\bbX^*$, i.e., the dual space of $\bbX$. % and %$h:\bbX\to\bbR\cup\{+\infty\}$ is a proper, closed and convex (but possibly non-smooth) function  such that 
Finally, $\emptyset\ne \calX\subseteq\bbX$ is a (bounded) polytope such that $\rvA(\calX)\subseteq\calK$ and $\rvA(\calX)\cap\inter\calK\ne \emptyset$, where $\rvA(\calX)$ %:= \{\rvA x:x\in\calX\}$ 
denotes the image of $\calX$ under $\rvA$.  {Under these assumptions, it follows that $\poi$ has at least one optimal solution, and we denote the set of optimal solutions by $\calX^*$. %Therefore, $F^*=F(x^*)$, for any $x^*\in \calX^*$. }
%The problem~\eqref{eq:poi} has many applications, including D-optimal design, Poisson linear inverse problem, the PET problem, and reformulated TV-regularized Poisson image denoising. For detailed description of these applications, we refer readers to~\cite{Zhao_23}. 

The FW method was first proposed in the seminal paper of Frank and Wolfe in 1956~\cite{Frank_56}, and received a fair amount of research interest until 1980s (see e.g.,~\cite{Lev_66,dem1967minimization,Canon_68,Dunn1980}). During the past decade, this method has attracted a resurgent and significant amount of research efforts (see e.g.,~\cite{Jaggi_13,Harcha_15,Freund_16,Freund_17,Nest_18,Ghad_19}), for at least two reasons. First, compared to gradient-projection-type methods that involve solving projection sub-problems onto $\calX$ (see e.g,~\cite{Nest_13}), it only requires solving linear-optimization  sub-problems over $\calX$, which may 
incur lower computational cost for certain types of  $\calX$ (e.g., polytopes or spectrahedra). Second, under certain scenarios, it can produce a sparse or low-rank solution, which is desired in many statistical and imaging applications. Specifically, if $\calX$ is a unit simplex (resp.\ spectrahedron), % = \Delta_m:= \{x\in\bbR^m:\sum_{i=1}^m x_i=1,x\ge 0\}$ (resp.), %namely the $(m-1)$-dimensional unit simplex, %and starting point is a
the number of non-zeros (resp.\ rank) of the iterates only increases by at most one in each iteration of the FW method. %Similarly, if $\calX$ is a spectrahedron, each iteration of the FW method only increases the rank of the (matrix) iterates by at most one. 

Despite its advantages, the FW method is also known to exhibit ``zig-zag'' behavior and hence suffer from slow convergence when the optimal solution lies on the relative boundary of $\calX$ (see e.g., Wolfe~\cite[Fig.\ 2]{Wolfe_70}). To remedy this drawback, Wolfe~\cite[Section 8]{Wolfe_70} proposed his ``away-step'' strategy, resulting in the away-step FW method. Both the algorithm and its convergence analysis were later corrected and refined by Gu\'elat and Marcotte~\cite{Guelat_86}. In particular, they showed that
when the objective function $F$ is both smooth and strongly convex on $\calX$, and strict complementarity holds for~\eqref{eq:poi}, 
the away-step Frank-Wolfe method identifies the minimal optimal face of $\calX$ and exhibits local linear convergence. 
Here we call $F$ smooth on $\calX$ if it is differentiable on some open set containing $\calX$ and its gradient, denoted by $\nabla F$, %of $F$ 
is Lipschitz   on $\calX$, namely for some $L\ge 0$, we have $\normt{\nabla F(x) - \nabla F(x')}_*\le L\normt{x-x'}$ for all $x,x'\in\calX$, where  $\normt{\cdot}_*$ denotes the dual norm of $\normt{\cdot}$. 
%By ``$F$ being smooth on $\calX$'', we specifically refer to the gradient of $F$ being Lipschitz   on $\calX$. %, namely  %w.r.t.\ certain norm $\normt{\cdot}$ and its dual 
More recently, several authors~\cite{Lacoste_15,Beck_17,Pena_19,Garber_20} provided global linear convergence analysis for the away-step FW method, when $F$ satisfies the smoothness and quadratic-growth conditions on $\calX$.

Let us note that all the algorithms and analyses mentioned above focus on the case where $F$ is smooth on $\calX$, which has become a rather standard setting for first-order methods (FOM) in general~\cite{Nest_04}. In this work, we instead focus on a very %fundamentally 
different problem class as described in~\eqref{eq:poi}, wherein $F$ may ``blow-up'' on (parts of) the relative boundary of $\calX$, and hence is not smooth (or %more generally, 
H\"older-smooth~\cite{Nest_15}) on $\calX$. 
Problems with %this type of 
such ``non-standard'' behaviors have recently attracted %certain amount of 
some research efforts.  %which we shall review in Section~\ref{sec:lit}. 
%This type of problems
In particular,  in the pivotal work~\cite{Dvu_20}, Dvurechensky et al.\ first proposed a novel FW method for  a similar problem class of~\eqref{eq:poi}, %in the sense that 
where $F$ is required to be a {non-degenerate} self-concordant function (but not necessarily a barrier nor logarithmically homogeneous). As an outgrowth of~\cite{Dvu_20}, Zhao and Freund~\cite{Zhao_23} %took this line of research into another direction, by 
proposed a generalized FW method for solving the problem %considered the case 
where $F$ is the sum of a (potentially) degenerate $\theta$-LHSCB and a ``simple'' closed convex function. 
%Indeed, 
Apart from this line of research, several different ``non-standard'' problem classes have been considered in recent works (e.g., Bauschke et al.~\cite{Bauschke_17}, Lu et al.~\cite{Lu_18} and Zhao~\cite{Zhao_22mg}), wherein principled FOMs were proposed to solve these problems --- for a detailed review, see Section~\ref{sec:lit}. 
%other authors have considered  and 

Before motivating our work, let us first formally review the definition of $\theta$-LHSCB. 
 %(see e.g.,). 

%for certain high-dimensional problems, if starting from a sparse (or low-rank) feasible point, the solution that it returns also has the sparse (or low-rank) property, provided that the number of iterations does not grow too large. 
%be  computationally advantageous for certain types of constraint set $\calX$ (for details, see e.g.,~\cite{Jaggi_13}). 

%\subsection{Logarithmically-homogeneous self-concordant barrier} 
\subsection{$\theta$-LHSCB}\label{sec:LHSCB}

%We let $f$ %that encompasses the $D$-optimal design problem, %arise in practice, 
% be a $\theta$-
%Let us formally review the definition of a logarithmically-homogeneous self-concordant barrier. %, which we now review.  
  Let $\calK\subseteq \bbY$ be a regular cone, i.e., $\calK$ is closed, convex, pointed (i.e., contains no straight line) and has nonempty interior (i.e., $\inter\calK\ne \emptyset$).  We say that $f$ is a $\theta$-logarithmically-homogeneous self-concordant barrier on $\calK$ for some $\theta \ge 1$, % and write $f \in {\sf LHB}_\theta(\calK)$, 
  if it satisfies the following properties:  
\begin{enumerate}[label=(P\arabic*),leftmargin=3.9em]
\item $f$ is convex and three-times continuously differentiable on $\inter \calK$,  \label{item:strict_convexity}
\item $\abst{D^3f(y)[u,u,u]}\le 2 \normt{u}_y^3$, \;\;\mbox{where} $\normt{u}_y:=\lranglet{\nabla^2 f(y)u}{u}^{1/2}$,  \quad $\forall\,y\in\inter\calK$, $\forall\,u\in\bbY$,\label{item:third_second_bounded}
\item $f(y_k)\to  +\infty$ for any $\{y_k\}_{k\ge 1}\subseteq\inter\calK$ such that $y_k\to y\in\bdry\calK$, and\label{item:boundary_growth}
\item $f(ty) = f(y) - \theta\ln (t)$ \  $\forall\,y\in\inter\calK$, $\forall\,t>0$ . \label{item:log_homogeneous}
\end{enumerate}
%where $H(x)$ is the Hessian of $f$ at $x$, see \\
(We refer interested readers to%Nesterov and Nemirovski
~\cite[Section 2.3.3]{Nest_94} and \cite[Section~2.3.5]{Renegar_01} for more details.)  
%We will assume in \eqref{eq:poi} that $f \in \calB_\theta(\calK)$ and that $\calX \subseteq \calK$.   
Properties~\ref{item:strict_convexity},~\ref{item:third_second_bounded} and~\ref{item:boundary_growth} correspond to $f$ being a (standard and strongly) self-concordant function with $\dom f = \inter\calK$  %with $\dom f = \inter\calK$ 
(cf.~\cite[Remark~2.1.1]{Nest_94}). In addition, since $\calK$ is pointed, we know that $f$ is non-degenerate, namely $\nabla^2 f(x)$  is positive definite for all $x\in\inter \calK$ (cf.~\cite[Theorem~4.1.3]{Nest_04}). We denote the class of  (standard, strongly and non-degenerate) self-concordant functions on $\inter\calK$ by ${\sf SC}(\calK)$.  If $f$ additionally satisfies property~\ref{item:log_homogeneous}, then it  is a $\theta$-LHSCB on $\calK$. 
%logarithmically-homogeneous barrier function on $\calK$, and we write $f\in {\sf LHB}_\theta(\calK)$. %(Therefore, ${\sf LHB}_\theta(\calK)\subseteq {\sf SC}(\calK)$.) 
Here $\theta$ is called the {\em complexity parameter} of $f$ (cf.\ Renegar~\cite{Renegar_01}). Finally, we let $f(x)=+\infty$ for $x\not\in\inter\calK$. %following the convention in convex analysis, %note that for $x\not\in\inter\calK $, . 

%From the definition above, 
Let us notice that the class of $\theta$-LHSCB functions is {\em fundamentally different} from the class of smooth functions. %in the sense that 
Specifically, its definition does not involve any pre-specified norm --- instead, as~\ref{item:third_second_bounded} indicates, %suggests, 
the norm is only defined locally via the Hessian of $f$. %$\nabla^2 f(y)$, and it changes as $y$ changes. 
Furthermore, at any $y\in\inter\calK$, $f$ only has an upper curvature bound in a ``neighborhood'' of $y$ (see Lemma~\ref{lem:ub_f} for details),
%the precise statement),  % 
which is %drastically 
very different from the {\em global} quadratic upper bound of any smooth function. 
%the Hessian of $f$. 
%In addition, %

\subsection{Motivating Example: D-Optimal Design} \label{sec:D-opt} % (and its variants)}

Given $m$ points $\{a_i\}_{i=1}^m \subseteq \mathbb{R}^n$ whose linear span is $\mathbb{R}^n$, the %(continuous) 
D-optimal design problem reads
\begin{align}
\min \;  F(x):= -\ln\det \big(\textstyle\sum_{i=1}^m x_i a_ia_i^T\big) \quad \st \;\; x\in\Delta_m, \label{eq:Dopt} \tag{DOpt}
\end{align}
where  $m\ge n+1$ and 
$$\textstyle\Delta_m:= \{x\in\bbR^m:\sum_{i=1}^m x_i= 1,\; x\ge 0\}$$ denotes the $(m-1)$-dimensional (standard) unit simplex. % and $e:=(1,\ldots,1)$. % denotes the vector of all ones. 
In statistics, this problem is the continuous relaxation of the well-known (discrete) D-optimal design problem:
\begin{align}
\min \;  -\ln\det \big(\textstyle\sum_{i=1}^m x_i a_ia_i^T\big) \quad \st \;\; x_i\in \bbZ_+, \;\;\sum_{i=1}^m x_i = k, \label{eq:Dopt_d}
\end{align}
where $\bbZ_+$ denotes the nonnegative integers, $k\in \bbZ_+$ and $k\ge n$. 
%The problem~\eqref{eq:Dopt_d} corresponds to 
In words, we wish to select (with repetition) $k$ design points from $\{a_i\}_{i=1}^m $  to maximize the determinant of the resulting Fisher information matrix $\sum_{i=1}^m x_i a_ia_i^T$~\cite{Fedorov_72}.  In computational geometry, $D$-optimal design arises as a Lagrangian dual problem of the minimum volume covering ellipsoid (MVCE) problem, which dates back at least 70 years to \cite{john} (see also Todd~\cite{toddminimum}).\\[-2ex]

\noindent {\bf Budget-constrained extension}.\ In certain situations, one also faces some budget constraints in designing the experiments, which leads to the following budget-constrained variant of~\eqref{eq:Dopt}:
\begin{align}
\min -\ln\det \big(\textstyle\sum_{i=1}^m x_i a_ia_i^T\big) \quad \st \;\; x\in\Delta_m, \;\; c_l^\top x \le C_l, \;\;\forall\,l\in[L].\label{eq:Dopt2}
\end{align}
In~\eqref{eq:Dopt2}, each constraint $c_l^\top x \le C_l$ corresponds to the limitation of one resource (e.g., time, money and labor), and we have $c_l\ge 0$ and $C_l> 0$, for all $l\in [L]$. Of course, we assume that~\eqref{eq:Dopt2} has at least one feasible solution. 
Here we see that both~\eqref{eq:Dopt} and~\eqref{eq:Dopt2} is an instance of \eqref{eq:poi} with $\calK=\bbS_+^n$, $f = -\ln\det(\cdot)$, $\theta = n$, $\rvA: x \mapsto \textstyle\sum_{i=1}^m x_i a_ia_i^T$ and $c=0$. In addition,  the constraint set $\calX\subseteq\Delta_m$ is a nonempty polytope that satisfies $\rvA(\calX) \subseteq \bbS_+^n$.\\[-2ex]

\noindent {\bf The WA-TY method, and an open question}. In 1973, Atwood~\cite{Atwood_73} proposed a modification of Fedorov's algorithm~\cite[Section~2.5]{Fedorov_72} for solving~\eqref{eq:Dopt}, which %has the same structure as 
structurally can be regarded as the FW method with exact line-search. %Atwood's modification 
It turns out that Atwood's modification coincides with Wolfe's away-step strategy~\cite{Wolfe_70}, and was later rediscovered by Todd and Y{\i}ld{\i}r{\i}m~\cite{Todd_07}. Following the convention of~\cite{Ahi_08}, we call it the WA-TY method. In practice, the WA-TY method has remarkable numerical performance  %for solving~\eqref{eq:Dopt}, 
--- in particular, it leads to significant speed-up over Fedorov's algorithm (see e.g.,~\cite{Atwood_73,Ahi_08}). Indeed, Ahipasaoglu et al.~\cite[Figure~1]{Ahi_08} observed numerically that the WA-TY method has global linear convergence (in terms of FW gap), % the sequence of FW gaps), 
however, they (only) managed to prove the local linear convergence. Indeed, establishing the global linear convergence of the WA-TY method remained an open question. 

%In fact, 
In~\cite[Section~3]{Ahi_08}, two structural difficulties of~\eqref{eq:Dopt} have been identified:
%, namely, 
i)  %to prove the global linear convergence is precisely that 
$F$ is not smooth on $\Delta_m$ (and ``blows-up'' on certain faces of $\Delta_m$) and ii) $F$ is degenerate on the feasible region (i.e., $\nabla^2 F(x)$ is non-invertible for any feasible $x$). 
An important motivation of the current work is to resolve the aforementioned open question, by overcoming the two difficulties above. As a key component of our approach, we identify the essential structure of~\eqref{eq:Dopt} that drives the global linear convergence of the WA-TY method, that is, $-\ln\det(\cdot)$ being a $n$-LHSCB. This property is fundamentally different from its smoothness counterpart, which is commonly assumed in the literature of the FW-type methods (and %first-order methods 
FOM in general). As a result, we are able to 
%and This enables us to 
consider a general class of problems as described in~\eqref{eq:poi}, which not only subsumes~\eqref{eq:Dopt}, but also  includes many other applications, among which we detail two in the following. 
% of the D-optimal design problem~\eqref{eq:Dopt} that 

\subsection{Other Applications} \label{sec:two_other}

\subsubsection{Inference of multivariate Hawkes processes (MHP)~\cite{Hawkes_71,Zhou_13}}

%\noindent {\bf }. %An $m$-dimensional multivariate Hawkes process can be described as follows. 
We adopt the description of an $m$-dimensional multivariate Hawkes process as in~\cite{Zhao_22w}. 
Fix a time interval $[0,t)$ for some $t>0$, and let $\calD:=\{(t_i,h_i)\}_{i=1}^n$ be the arrival points in $[0,t)$, where for each $i\in[n]$, $t_i$ denotes the arrival time and $h_i\in[m]$ denotes the index of the dimension along which the point arrives. We assume that arrival points occur along every dimension $k\in[m]$, which happens with high probability for $t$ large enough. 
For each dimension $k\in[m]$, the conditional density function is given by
\begin{equation}
\lambda_k(t):= \mu_k + \textstyle{\sum}_{i=1}^n\; a_{h_i,k}\exp(-(t-t_i))\bbone_{>0}(t-t_i) , \quad\forall\, t>0, \label{eq:lambda_k}                                                                      
\end{equation}
where $\mu_k\ge 0$ is the base intensity at dimension $k$ and %$\{a_{h_i,k}\}_{i=1}^n$ are the
$a_{h_i,k}\ge 0$ is the mutual-excitation coefficient from dimension $h_i$ to dimension $k$. In addition, $\bbone_{>0}$ denotes the step function, i.e., $\bbone_{>0}(a) = 1$ for $a>0$ and $0$ for $a\le 0$.  %, and $\zeta:\bbR_+\to\bbR_+$ is the kernel function. 
We are interested in inferring the base intensity vector $\bar \mu:= (\mu_k)_{k\in[m]}$ and the infectivity matrix $A:= (a_{k,l})_{k,l\in[m]}$ via maximum-likelihood estimation (MLE), based on the arrival points $\calD$ in $\calT$. Note that  the log-likelihood function  reads (see e.g.,~\cite{Zhou_13})
\begin{equation}
%\begin{split}
L(\bar\mu,A) := \textstyle\sum_{i=1}^n \ln \lambda_{h_i}(t_i) - \sum_{k=1}^m \int_{0}^t\lambda_k(s)\rmd s. 
%\\
%& = \textstyle\sum_{i=1}^n \ln \big(\mu_{h_i} + \sum_{j:t_j<t_i} a_{h_j,h_i}\, \zeta(-(t_i-t_j))\big)\\
%& \quad \textstyle - \sum_{k=1}^m \big(\mu_k t + \sum_{i=1}^n a_{h_i,k}\, \int_{0}^t \zeta(-(s-t_i))\rmd s\big).
%\end{split}
\end{equation}
Now, let us partition the arrival points $\calD$ along $m$ dimensions. Specifically,  for each $k\in[m]$, define the index set $\calH_k(t):= \{i\in[n]:t_i<t,h_i=k\}\ne \emptyset$, and then $\{(t_i,h_i)\}_{i\in \calH_k(t)}$ denote the arrival points along dimension $k$. Using this notation, we can write $L$ in a dimension-separable form:
\begin{equation}
\textstyle L(\bar\mu,A) = \sum_{k=1}^m L_k(\mu_k,a_k), \quad\mbox{where}\quad L_k(\mu_k,a_k) := \sum_{i\in\calH_k(t)} \ln \lambda_{h_i}(t_i) - \int_{0}^t\lambda_k(s)\rmd s,  
%  \big(\mu_{k} + \sum_{l=1}^m  a_{l,k}\,\sum_{j\in\calH_l(t_i)} \zeta(-(t_i-t_j))\big)
\end{equation}
and $a_k^\top$ denotes the $k$-th row of $A$. This suggests that the inference of $\bar \mu$ and $A$ can be carried out %independently {in parallel} 
{\em separately} along the $m$ dimensions, and hence let us focus on a single dimension $k\in[m]$. Using the definition of $\lambda_k$ in~\eqref{eq:lambda_k}, we can write $L_k$ in the following form:
%\begin{align*}
%
%\end{align*}
%
%so that $\{\calH_k(t)\}_{k\in[m]}$ form a partition of $[n]$. We assume that $\calH_k(t)\ne\emptyset$ for all $k\in[m]$. (This happens with high probability if $t$ is large enough.) We then observe that the log-likelihood function $L(\cdot,\cdot)$ is separable across the $m$ dimensions, namely $L(\mu,A):= \sum_{k=1}^m L_k(\mu_k,a_k)$, where $a_k^\top$ denotes the $k$-th row of $A$ and 
%\begin{equation}
%\begin{split}
%L_k(\mu_k,a_k) & = \textstyle\sum_{i\in\calH_k(t)} \ln \big(\mu_{k} + \sum_{l=1}^m  a_{l,k}\,\sum_{j\in\calH_l(t_i)} \zeta(-(t_i-t_j))\big)\\
%& \quad \textstyle - \big(\mu_k t + \sum_{l=1}^m a_{l,k}\, \sum_{i\in\calH_l(t)}\int_{0}^t \zeta(-(s-t_i))\rmd s\big).
%\end{split}
%\end{equation}
%In particular, if we let $\zeta(t) = \exp(-t)$ for $t>0$ and $0$ for $t\le 0$, then we have
\begin{equation}
\begin{split}
L_k(\mu_k,a_k)& = \textstyle\sum_{i\in\calH_k(t)} \ln \big(\mu_{k} + \sum_{l=1}^m  a_{l,k}\,\sum_{j\in\calH_l(t_i)}   \exp(-(t_i-t_j))\big)\\
& \quad \textstyle - \big( t \mu_k + \sum_{l=1}^m a_{l,k}\, \sum_{i\in\calH_l(t)}(1-\exp(-(t-t_i)))\big)
\\
&= \textstyle \sum_{i\in\calH_k(t)} \ln \big(\mu_{k} + \sum_{l=1}^m  a_{l,k} \,(\barw_{i})_{l} \big) - (t \mu_k + \sum_{l=1}^m a_{l,k}\, v_l),\\
\end{split}
\end{equation}
where for all $i\in\calH_k(t)$ and $l\in[m]$, 
\begin{equation}
(\barw_{i})_l:= \textstyle\sum_{j\in\calH_l(t_i)}   \exp(-(t_i-t_j))\ge 0  \quad\mbox{and}\quad  v_l:= \sum_{i\in\calH_l(t)}1-\exp(-(t-t_i)) >0, 
\end{equation}
where we use the convention that empty sum is zero. By denoting $\barw_i:= (\barw_{i,l})_{l\in[m]}$ and $v:= (v_{l})_{l\in[m]}$,
%$\psi_l:= \sum_{i\in\calH_l(t)}(1-\exp(-(t-t_i)))$, . 
%Note that since $\calH_l(t)\ne \emptyset$, we have $v_l> 0$, for all $l\in[m]$. 
 the MLE problem along the $k$-th dimension reads
\begin{equation}
{\min}_{\mu_k\ge 0,\, a_k\ge 0} \; \textstyle -\sum_{i\in\calH_k(t)} \ln (\mu_{k} +   \barw_{i}^\top a_{k} ) +(t \mu_k +  v^\top a_{k} ). \label{eq:MLE}
\end{equation}
Sometimes, to promote the sparsity of the coefficients $\{a_{l,k}\}_{l=1}^m$, we add an $\ell_1$-regularizer to~\eqref{eq:MLE}:
\begin{equation}
\textstyle{\min}_{\mu_k\ge 0,\, a_k\ge 0} \;  -\sum_{i\in\calH_k(t)} \ln \big(\mu_{k} +   \barw_{i}^\top a_{k}  \big) +(t \mu_k +  v^\top a_{k}) + \lambda \normt{a_{k}}_1, \label{eq:L1_MLE}
\end{equation}
where $\lambda\ge 0$ and $\normt{a_{k}}_1:= \sum_{l=1}^m \abst{(a_{k})_l} = \sum_{l=1}^m {(a_{k})_l}$, since $a_k\ge 0$. Using the standard re-scaling technique (see e.g.,~\cite[Section~2.2]{Zhao_22w}), we can write~\eqref{eq:L1_MLE} equivalently in the following form:
\begin{equation}
\begin{split}
{\min}_{(\mu,a)\in\bbR\times \bbR^m} \quad &\textstyle F(\mu,a):=-\sum_{i\in\calH_k(t)} \ln \big(\mu/t +   w_i^\top a \big)\\
 \st \quad &(\mu,a)\in \Delta_{m+1}:=\{(\mu,a)\in\bbR\times \bbR^m: \mu\ge 0,\;\;a\ge 0,\;\;\mu + e^\top a = 1\}. \label{eq:L1_MLE_rewrite}
\end{split} \tag{MHP}
\end{equation}
In~\eqref{eq:L1_MLE_rewrite}, $e:= (1,\ldots,1)$ and $w_i:= \barw_i \oslash (v+\lambda e)$ for $i\in\calH_k(t)$, where $\oslash$ denotes entrywise division. 
%In addition, note that we can write the constraint in~\eqref{eq:L1_MLE_rewrite} as $(\mu,a)\in\Delta_{m+1}$. 
Note that the problems in~\eqref{eq:L1_MLE} and~\eqref{eq:L1_MLE_rewrite} are equivalent in the following sense: 
$(\mu^*,a^*)$ is an optimal solution of~\eqref{eq:L1_MLE_rewrite} if and only if $(\mu^*/t,a^* \oslash(v+\lambda e) )$ is an optimal solution of~\eqref{eq:L1_MLE}.   % and we know that $v+\lambda e>0$. 
% = e^\top a_{k}$ and $e:= (1,\ldots,1)$ . 
%Note that since $a_k\ge 0$, we have $\normt{a_{k}}_1=e^\top a_{k}$, where $e:= (1,\ldots,1)$. %The formulation in~\eqref{eq:L1_MLE} can also be interpreted as maximum a-posteriori estimation of $a_k$ with independent exponential prior on the elements $\{a_{l,k}\}_{l=1}^m$. 

Indeed, the problem~\eqref{eq:L1_MLE_rewrite} can be written in the following simple abstract form:
\begin{equation}
{\min}_{z} \quad \textstyle -\sum_{i=1}^m \ln \big(a_i^\top z \big)\quad  \st \quad z\in\Delta_n, \label{eq:PET}
\end{equation}
which is an instance of \eqref{eq:poi} with $\calK=\bbR^m_+$, $f(y) = -\sum_{i=1}^m \ln(y_i)$, $\theta = m$, $\rvA: x \mapsto A^\top x$ for the $n\times m$ data matrix $A:=[a_1\;\ldots\; a_m]$ and $c=0$. In addition,  the constraint set $\calX = \Delta_n$, and since $A$ is (entrywise) nonnegative, we have $\rvA(\calX)\subseteq \bbR^m_+$. Somewhat surprisingly, the problem~\eqref{eq:PET} also appears in applications across different fields, including positron emission tomography~\cite{Vardi_93}, sparse Poisson inverse problem~\cite[Section~5]{Bauschke_17} and (offline) log-optimal investment problems~\cite{Cover_84}. Note that depending on different applications, %depending on the specific application, 
the data matrix $A$ may have different sparsity structures. %may be different for these applications.) 
%For detailed descriptions, we refer readers to~\cite[Section~1.1]{Zhao_23}. 
\\[-2ex]

\subsubsection{Poisson image de-blurring with total-variation (TV) regularization~\cite{Harmany_12,Chambolle_18}} % see also \cite{Dvu_20}, \cite{Chambolle_11}. \textcolor{red}{[Renbo, let us discuss these references for this topic]}. In image processing 
%The following description is adopted from~\cite[Section~1.1]{Zhao_23}. 
Let the $m\times n$ matrix $X$ represent an image, such that each entry $X_{ij}\in\{0,1,\ldots,M\}$ represents the intensity of the pixel at location $(i,j)\in[m]\times[n]$. %, where $M:= 2^b-1$.  %for $i\in[m]$ and $j\in[n]$, 
%and $X_{ij}\in\{0,1,\ldots,M\}$, where $M:=2^b-1$ for $b$-bit images. % and $b\ge 1$ is an integer. %\in\bbN$, the set of natural numbers. %, , where $M$. % for $k$-bit images))
In many applications, ranging from microscopy to astronomy, we observe a blurred noisy image $Y\in\bbR_+^{m\times n}$, and we wish to estimate the true image from $Y$. Following the description as in~\cite[Section~1.1]{Zhao_23}, this leads to the following optimization problem: 
\begin{align}
{\min}_{x\in\bbR^N}&\;\; -\textstyle\sum_{l=1}^{N} y_l\ln(a_l^\top x) + (\sum_{l=1}^{N} a_l)^\top x + \lambda {\rm TV}(x)\quad \st\;\; 0\le  x\le Me \ , %\;\forall\,j\in[p]. 
\label{eq:deblurring_TV}
\end{align} %\red{[Renbo you had $\bar F(x)$ above, why the bar?]}
where $N:= mn$, $x\in\bbR^N$ and $y\in\bbR_+^N$ denote the vector representations of $X$ and $Y$, respectively, $a_l$  denotes the $l$-th column of an $N\times N$ nonnegative sparse (low-rank) matrix, and 
$${\rm TV}(x):= \textstyle\sum_{i=1}^m \sum_{j=1}^{n-1} \abs{X_{i,j} - X_{i,j+1}} + \textstyle\sum_{i=1}^{m-1} \sum_{j=1}^{n} \abs{X_{i,j} - X_{i+1,j}}.$$ % (recall that $X$ is the matrix representation of $x$). 
Clearly, by properly defining a (sparse) matrix $D\in\{0,\pm 1\}^{E\times N}$, where $E:= (n-1)m+(m-1)n$, we can write %${\rm TV}(x) = \normt{Dx}_1$, and consequently, represent ${\rm TV}(x)$ as the optimal value of the following LP:
\begin{align}
%\begin{split}
{\rm TV}(x)= \normt{Dx}_1 = {\min}_{r\in\bbR^E} %{\min}_{x\in\bbR^N,\;r\in\bbR^m} 
\;\; e^\top r\;\; \st \;\;r\ge D x, \;\; r\ge -D x, \;\;r\le Me. \label{eq:TV_LP} 
\end{align}
(Note that the constraint $r\le Me$ is in fact redundant, however, it ensures the boundedness of $r$.)
%$r\in\bbR^E$ for $E:= (n-1)m+(m-1)n$ and $B\in\bbR^{N\times E}$ is the node-arc incidence matrix of some directed network (see e.g.,~\cite[Section 2]{Goldfarb_09}). 
Using~\eqref{eq:TV_LP}, we can re-write~\eqref{eq:deblurring_TV} as 
\begin{align}
\begin{split}
{\min}_{x,r}%{x\in\bbR^N,\;r\in\bbR^m} 
\;\; &-\textstyle\sum_{l=1}^{N} y_l\ln(a_l^\top x) + \textstyle(\sum_{l=1}^{N} a_l)^\top x + \lambda e^\top r \quad \;\;\\
\st\;\; &0\le  x\le M e, \;\;r\ge D x, \;\; r\ge -D x, \;\;r\le Me. 
\end{split}\label{eq:box_TV2}
\end{align}
Note that this is an instance of~\eqref{eq:poi}, by taking $\calK=\bbR^N_+$, $f(z) = -\sum_{i=1}^N \ln(z_i)$, $\theta = N$, $\rvA: (x,r) \mapsto Ax$ and $c:(x,r) \mapsto \textstyle(\sum_{l=1}^{N} a_l)^\top x + \lambda e^\top r$. %$c=[e^\top A\;\;  \lambda e^\top]^\top$. 
In addition,  the constraint set $\calX:=\{x\in\bbR^N,r\in\bbR^E:0\le x\le M e, \;r\ge D x, \;r\ge -D x, \;r\le Me\}$ is a nonempty polytope and we have $\rvA(\calX)\subseteq \bbR^N_+$.
%In fact, it is well-known that ${\rm TV}(\cdot)$ can be represented as follows (see e.g.,~\cite[Section 6.2]{Harcha_15}): 

\subsection{Related Work}\label{sec:lit}

%As mentioned in Section~\ref{sec:LHSCB}, 
Since the problem class~\eqref{eq:poi} greatly differs from the ``standard'' setting %one commonly considered 
in the FOM literature, where $f$ is assumed to be smooth on $\calX$, it has not been extensively studied in the context of FOM.  %As a result, it has been rarely studied %That said, 
%In fact, 
To the best of our knowledge, Dvurechensky et al.~\cite{Dvu_20} pioneered the development of FOM for~\eqref{eq:poi} --- % Specifically, %, wherein 
they proposed a novel FW method for a similar problem class of~\eqref{eq:poi}, %in the sense that 
where $F$ is required to be a {non-degenerate} self-concordant function (but not necessarily a barrier nor logarithmically homogeneous). 
%Nevertheless, two fairly recent works have proposed some new variants of the FW-type method to solve~\eqref{eq:poi}. Specifically, 
They showed their algorithm converges at rate $O(C/k)$ in terms of objective gap for some positive constant $C$.  
Subsequently, the authors of~\cite{Zhao_23} took this line of research in another direction, by considering the case where $F$ is the sum of a (potentially) degenerate $\theta$-LHSCB and a ``simple'' closed convex function, and proposed a generalized FW method for solving this problem. % more general than~\eqref{eq:poi}, and 
When specialized to~\eqref{eq:poi}, they showed that the number of iterations to obtain both an $\varepsilon$-objective gap and an $\varepsilon$-FW gap is (essentially) of order  $O((\theta+B)^2/\varepsilon)$, where $B:= \max_{x,x'\in\calX}\,\ipt{c}{x-x'}$ denotes the variation of $\ipt{c}{\cdot}$ on $\calX$. In particular, on the D-optimal design problem, this  result (essentially) recovers the complexity result derived by Khachiyan~\cite{Khachiyan_96}. 
Meanwhile, Dvurechensky et al.~\cite{Dvu_23} extended their FW method in~\cite{Dvu_20}  to the case where $F$ is a  non-degenerate {\em generalized} self-concordant function~\cite{Sun_18}.  In the same work, they also proposed an away-step extension of their FW method when $\calX$ is polyhedral, and showed that under certain assumptions, their away-step FW method converges globally linearly %has global linear convergence 
in terms of the objective gap (see also~\cite{Card_21} for  a similar away-step FW method with backtracking line-search). 
%proposing a FW method for minimizing the  $F$ over a convex compact set $\calX\ne \emptyset$, and showed a convergence rate of $O(1/k)$. 

Motivated by establishing the global linear convergence of the WA-TY method on the D-optimal design problem (cf.\ Section~\ref{sec:D-opt}) as well as the two other applications in Section~\ref{sec:two_other},
%inference of MHP and TV-regularized Poisson image de-blurring), %, 
we were independently developing and analyzing our away-step FW method for solving~\eqref{eq:poi}. In %the final stage of completing 
finalizing the current manuscript, we became aware of the interesting results in~\cite{Dvu_23}. However, these results notably  % (and~\cite{Card_21}).While the results in~\cite{Dvu_23} are very interesting, it is unfortunate that they 
do not apply to any of the aforementioned applications, namely  %problem classes outlined herein, namely they cannot be applied to 
D-optimal design, inference of MHP, and TV-regularized Poisson image de-blurring.  The main reason is the requirement for $F$ to be {\em non-degenerate} (i.e., $\nabla^2 F(x)$ is invertible for all feasible $x$), which does not hold for any of the three applications. %problem classes of interest herein.
In fact, the degeneracy of $F$ was identified in  Ahipasaoglu et al.~\cite{Ahi_08} to be a key obstacle in  proving the  global linear convergence of the WA-TY method on D-optimal design. %~\eqref{eq:Dopt}. 
In addition, the away-step FW method in~\cite{Dvu_23} (implicitly) requires that at least one vertex of $\calX$ is  feasible for~\eqref{eq:poi}. While this holds for several applications, unfortunately, it does not hold for D-optimal design. 
%~\eqref{eq:Dopt}. %$\dom F \cap\Ext(\calX)\ne \emptyset$, where $\dom F:= \{x\in\bbX:F(x)<+\infty\}$  and $\Ext(\calX)$ denotes the vertices of $\calX$. 

% Upon careful examination, we found that these results do not apply to the D-optimal design problem and the other two applications in Section~\ref{sec:two_other}, for two reasons. First, these results assume that $\dom F \cap\Ext(\calX)\ne \emptyset$, where $\dom F:= \{x\in\bbX:F(x)<+\infty\}$  and $\Ext(\calX)$ denotes the vertices of $\calX$. Unfortunately, this assumption is not satisfied by~\eqref{eq:Dopt}. Second, and more importantly, the non-degeneracy assumption on $F$ does not hold on any aforementioned application. % mentioned previously. 
%In fact, as identified by Ahipasaoglu et al.~\cite{Ahi_08}, the degeneracy of $F$ is one of the major difficulties for proving the global linear convergence of the WA-TY method on~\eqref{eq:Dopt}, and %it does not admit an easy solution. %
% overcoming this difficulty requires considerable efforts (and indeed, more than we expect). 
%Fortunately, we manage to address these two issues in our algorithm and its analysis, and therefore resolve the open problem in~\cite{Ahi_08}. %and so 

In addition, some other aspects of our analysis are highlighted below. %compared to the re
%we also wish to highlight . % that we believe are valuable. 
First, the global linear rate of the objective gap that we obtain is both {\em affine-invariant} and {\em norm-independent}.\footnote{By ``norm-independent'', we mean the linear rate does not depend on any pre-specified norm, e.g., the $\ell_2$-norm.} 
As the away-step FW method is well-known to be affine-invariant and norm-independent (see e.g.,~\cite{Lacoste_15}), it is both natural and desirable for the convergence analysis  to have such properties as well %be affine-invariant  
(see e.g., Pe\~{n}a~\cite{Pena_21}). Second, we also establish the global linear convergence of the sequence of FW gaps generated by our away-step FW method. Compared to the objective gap, which is typically unknown in practice, the FW gap is easily computable and serves as an upper bound of the objective gap. In fact, it is the FW gap that was observed to converge globally linearly on the D-optimal design problem in Ahipasaoglu et al.~\cite[Figure~1]{Ahi_08}, and our result provides the first theoretical explanation for this. Third, we show that the iterates generated by our method (with exact line-search) will land on and remain in a face of $\calX$ within a bounded number of iterations, and the linear convergence rates (in terms of both the objective gap and the Frank-Wolfe gap) may become faster. Indeed, this result explains our numerical observations very well (see Section~\ref{sec:experiments}).

%In addition, they proposed an away-step extension for polyhedral $\calX$,   by assuming that   Moreover, note that none of the applications in Sections~\ref{sec:D-opt} and~\ref{sec:two_other} satisfies the {\em non-degeneracy} assumption of $F$. Moreover, this assumption may not be satisfied by the general problem in~\eqref{eq:poi}, since importantly, the linear operator $\rvA$ may not be injective. % in these applications). 
%As a result, the algorithms and analyses in~\cite{Dvu_23} may not be applicable to the problem class in~\eqref{eq:poi}.  
%%the problem considered in this paper.  % prevent their away-step extension (and the associated analysis)  to 

%At a broader level, 
To conclude this section, let us also briefly mention some other principled first-order methods that have been developed recently for convex differentiable optimization without smoothness property. Specifically, Bauschke et al.~\cite{Bauschke_17} and Lu et al.~\cite{Lu_18} considered the setting where the objective function $F$ is ``relatively smooth'' w.r.t.\ some reference function, and proposed the relatively-smooth gradient method.  %to solve it. 
In addition,  Zhao~\cite{Zhao_22mg} proposed a generalized multiplicative gradient method for minimizing a ``gradient log-convex'' function over a slice of symmetric cone. Note that both  methods converge at rate $O(1/k)$ in terms of objective gap, and can be used to solve some aforementioned problems, including D-optimal design in~\eqref{eq:Dopt} and inference of MHP in~\eqref{eq:L1_MLE_rewrite}. % some problems in Sections~\ref{sec:D-opt} and~\ref{sec:two_other}, including D-optimal design and . 
%for solving  lines of works that aim to develop 

% in the se

\subsection{Main Contributions}

%Our main contributions are summarized as follows.

First, we analyze an away-step FW method for solving~\eqref{eq:poi}, and show the global linear convergence of this method in terms of the objective gap. Our linear convergence rate %has the important property of being 
is {\em affine-invariant} and {\em norm-independent}. Indeed, it only depends on five affine-invariant quantities (one of which is a geometric constant in a similar spirit to the ``facial distance'' constant in Pe\~{n}a and Rodriguez \cite{Pena_19}), and does not depend on any {pre-specified} norm (e.g., the $\ell_2$ norm). Our analysis makes crucial and extensive use of the logarithmic homogeneity and self-concordant property of $f$, and certain results (e.g., Lemma~\ref{keyinequality}) may be of independent interest. 
 % that are all . 
%As the away-step FW method is well-known to be 
%affine-invariant and norm-independent (see e.g.,~\cite{Lacoste_15}), it is both natural and desirable for the convergence analysis  to have such properties as well %be affine-invariant  
%(see e.g.,~\cite{Pena_21}). 

Second, we show that the sequence of FW gaps generated by our method is also globally linearly convergent. The global linear rate of the FW gap is also affine-invariant and norm-independent, but is worse compared to that of the objective gap. In particular, when specialized to the D-optimal design problem, this provides theoretical justification of (part of) the empirical observations in Ahipasaoglu et al.~\cite[Figure~1]{Ahi_08}.

Third, we show that the iterates generated by our method (with exact line search) will land on and remain in a face of $\calX$ in a finite number of iterations. Since this face may have a (much) smaller dimension compared to that of $\calX$, the away-step FW method may have improved local linear convergence rates (in terms of both the objective gap and the Frank-Wolfe gap). In addition, if $\calX $ is the unit simplex, the returned solution may be sparse. % (i.e., have a small number of non-zeros). 
%the returned approximate solution 

Fourth, we conduct numerical experiments on the problems of D-optimal design and inference of MHP. Our results  show that, compared to other first-order methods, our away-step FW method is significantly more efficient in reducing the objective gap, and can produce a much sparser solution. In addition, the observations made from these %(empirical) 
results corroborate our theory quite well. 

%and FW gap 
%and can but also 
\section{Preliminaries on Self-Concordant Functions}

Let $f\in{\sf SC}(\calK)$, namely, $f$ is a (standard, strongly and non-degenerate) self-concordant function on $\inter\calK$ that satisfies properties~\ref{item:strict_convexity},~\ref{item:third_second_bounded} and~\ref{item:boundary_growth}. For any $y\in\inter\calK$, let us define the local norm of $u\in\bbY$ at $y$ as $\normt{u}_y:=\lranglet{\nabla^2 f(y)u}{u}^{1/2}$, the local norm of $v\in\bbY^*$  at $y$ as $\normt{v}_{y}^*:= \sup_{\normt{u}_y=1}\,\ipt{v}{u} = \lranglet{\nabla^2 f(y)^{-1}v}{v}^{1/2}$, and the (open) Dikin ellipsoid at $y$ with radius $r> 0$ as
\begin{equation}
\calB_y(y,r):= \{u\in\bbY: \normt{u-y}_y< r\}. 
\end{equation}
%For any $v\in\bbY^*$ (i.e., the dual space of $\bbY$), define its 
Note that if $f\in{\sf SC}(\calK)$, then $\calB_y(y,1)\subseteq\inter\calK$ for all $y\in\inter\calK$ (see e.g., \cite[Theorem~4.1.5]{Nest_04}). Next, we will review some curvature bounds of $f\in{\sf SC}(\calK)$, which will be crucial in developing and analyzing our algorithms. To that end, let us first define the uni-variate (closed and convex) function 
\begin{equation}
\omega(t):= t- \ln(1+t), \quad \forall\, t> -1,
\end{equation}
and its Fenchel conjugate 
\begin{equation}
\omega^*(t):={\sup}_{s\in\bbR}\; ts - \omega(s) = -t - \ln(1-t), \quad \forall\, t<1. 
\end{equation}
Incidentally, we have $\omega^*(t) = \omega(-t)$ for all  $t<1$. 

\begin{lemma}[Lower curvature bound of $f$ {\cite[Theorem 4.1.7]{Nest_04}}]\label{lem:lb_f}
 If $f\in{\sf SC}(\calK)$, then we have %for all $x,y\in\inter\calK$, we have
 \begin{equation}
 f(y')\ge f(y) + \lranglet{\nabla f(y)}{y'-y} + \omega(\normt{y'-y}_y), \quad \forall\, y',y\in\inter\calK.
 \end{equation}
% where $\omega(t):= t - \ln(1+t)$ for $t\ge 0$. 
\end{lemma}

\begin{lemma}[Upper curvature bound of $f$ {\cite[Theorem 4.1.8]{Nest_04}}]\label{lem:ub_f}
 If $f\in{\sf SC}(\calK)$, then we have %for all $x,y\in\inter\calK$, we have
 \begin{equation}
 f(y')\le f(y) + \lranglet{\nabla f(y)}{y'-y} + \omega^*(\normt{y'-y}_y), \quad \forall\, y\in\inter\calK,\;\; \forall\, y'\in\calB_y(y,1). 
 \end{equation}
% where $\omega(t):= t - \ln(1+t)$ for $t\ge 0$. 
\end{lemma}

\begin{lemma}[Lower bound of $\omega$] \label{lem:omega}
Fix any $\beta>0$. Then we have %and any $t\in(0,\beta]$, we have $\omega(t)\ge \mu_\beta t^2$, where 
\begin{alignat}{3}
\omega(t)&\ge \mu_\beta t^2,\quad &&\forall\,t\in[0,\beta],\quad &&\mbox{where }\;\; \mu_\beta:= \omega(\beta)/\beta^2, \label{eq:quad_omega}\\
\omega(t)&\ge \varrho_\beta t,\quad &&\forall\,t\ge \beta,\quad &&\mbox{where }\;\; \varrho_\beta:= \omega(\beta)/\beta. \label{eq:lin_omega}
\end{alignat}
\end{lemma}
\begin{proof}
See Appendix~\ref{app:proof_omega}. 
%This completes the proof. 
% Since $\xi(0)=0$ and for all $t >0$, $\xi'(t) = -t^2/(1+t)^2<0$, we see that for all $t >0$, $\xi(t)<0$ and hence $\zeta'(t)<0$. %This completes the proof. 
\end{proof}

\begin{lemma}[Upper bound of $\omega^*$] \label{lem:ub_omega*}
We have 
\begin{equation}
\omega^*(t)\le \bigg\{\bar\omega^*(t):=\frac{t^2}{2(1-t)}\bigg\}, \quad \forall\, t\in[0,1).  \label{eq:omega*}
\end{equation}
\end{lemma}
\begin{proof}\let\qed\relax
Note that $\omega^*(0)=0$ and hence for all $t\in[0,1)$, we have 
\begin{equation*}
\omega^*(t) = \int_{0}^t (\omega^*)'(s)\; \rmd s = \int_{0}^t \frac{s}{1-s}\; \rmd s\le \int_{0}^t \frac{s}{1-t}\; \rmd s = \frac{t^2}{2(1-t)}. \tag*{$\square$} %\eqno\qed %
\end{equation*}%\noqed
\end{proof}

\begin{lemma}[Hessian approximation~{\cite[Theorem~4.1.6]{Nest_04}}] \label{lem:approx_Hessian}
For any $y\in\inter\calK$ and any $y'\in\calB_y(y,1)$, we have 
\begin{equation}
(1-\normt{y'-y}_y)^{2}\,\nabla^2 f(y)\preceq \nabla^2 f(y') \preceq(1-\normt{y'-y}_y)^{-2}\,\nabla^2 f(y). 
\end{equation}
\end{lemma}

\begin{lemma}\label{lem:local_norm_grad_diff}
For any $y\in\inter\calK$ and any $y'\in\calB_y(y,1)$, we have
\begin{equation}
\normt{\nabla f(y') - \nabla f(y)}_{y}^*\le \frac{\normt{y'-y}_y}{1-\normt{y'-y}_y}. 
\end{equation}
\end{lemma}

\begin{proof}
See Appendix~\ref{app:proof_grad_diff}. 
\end{proof}

\section{Away-step Frank-Wolfe Method for Solving~\eqref{eq:poi}}

 %Note that we need not 

\begin{algorithm}[t!]
\caption{Away-step Frank-Wolfe Method for solving~\eqref{eq:poi}}\label{algo:AFW_SC}
\begin{algorithmic}
\State {\bf Input}: Starting point $x^0\in\calX\cap\dom F$, the set of active atoms at $x^0$, denoted by $\calS_0\subseteq \calV$, and the atom-weight vector $\beta^0\in\Delta_{\abst{\calV}}$ such that $\supp(\beta^0) :=\{a\in\calV:\beta^0_a\ne 0\}= \calS_0$.\\\vspace{-1ex} %$\{\beta_a^0\}_{a\in\calS_0}$, such that $x^0 = \sum_{a\in\calS_0} \beta_a^0 a$, where $\beta_a^0>0$ for all $a\in\calS_0$ and $\sum_{a\in\calS_0} \beta^0_a = 1$. % and $\beta_a^0=0$ for all $a\not\in\calS_0$. % , number of iterations $K$ %accuracy $\varepsilon>0$
\State {\bf At iteration $k\ge 0$}:
\begin{enumerate}[leftmargin = 6ex]
\item   Compute $y^k := \rvA x^k$  %, $\nabla f(y^k)$ 
and $\nabla F(x^k) = \rvA^*\nabla f(y^k)+c$. 
\item \label{item:LMO_a} Compute $v^k\in\argmin_{x\in\calV}\;  \lranglet{\nabla F(x^k)}{x}$, $s^k := \rvA v^k$, the FW step $d_\rmF^k := v^k - x^k$ and 
\begin{equation}
G_k := \lranglet{-\nabla F(x^k)}{d_\rmF^k} =  \lranglet{\nabla f(y^k)}{y^k - s^k}+ \lranglet{c}{x^k - v^k}. \label{eq:G_k}
\end{equation}
If $G_k=0$, then STOP; otherwise proceed to Step~\ref{item:away_atom_a}. 
\item \label{item:away_atom_a} 
If $\abst{\calS_k}>1$, compute $a^k \in \argmax_{x\in\calS_k} \lranglet{\nabla F(x^k)}{x}$, $w^k := \rvA a^k$, the away step $d_\rmA^k := x^k - a^k$ and 
\begin{equation}
\tilG_k := \lranglet{-\nabla F(x^k)}{d_\rmA^k} = \lranglet{\nabla f(y^k)}{w^k - y^k}+ \lranglet{c}{a^k - x^k}. \label{eq:tilG_k}
\end{equation}
\item \label{item:choose_direction} Choose between the following two cases: 
\begin{enumerate}[label=(\alph*)]
\item If $\abst{\calS_k}=1$ or $G_k > \tilG_k$, let $d^k := d_\rmF^k$ and  $\baralpha_k := 1$.
\item Otherwise, let $d^k:= d_\rmA^k$ and $\baralpha_k := \frac{\beta^k_{a^k}}{1-\beta^k_{a^k}}$. 
\end{enumerate} 
\item \label{item:step_size_SC_a} Choose $\alpha_k\in(0,\baralpha_k]$ in one of the following two ways:
\begin{enumerate}[label=(\alph*)]
\item Exact line-search: \label{item:exact_LS} $\alpha_k\in \argmin_{\alpha_k\in(0,\baralpha_k]} F(x^k + \alpha d^k)$.
\item \label{item:adapt_step} 
Adaptive stepsize: Compute $r_k := -\ipt{\nabla F(x^k)}{d^k}$ and $D_k := \normt{\rvA d^k}_{y^k}$. If $D_k = 0$, then $\alpha_k := \baralpha_k$; otherwise, 
\begin{equation}
\alpha_k:= \min\{b_k,\baralpha_k\},\quad \mbox{where}\;\; b_k:=\frac{r_k}{D_k(r_k+D_k)}. \label{eq:def_alpha_k}
\end{equation}
\end{enumerate} 
\item Update $x^{k+1}:= x^k + \alpha_k d^k$, $\calS_{k+1}\subseteq \calV$ (i.e., the set of active atoms at $x^{k+1}$) and the atom weights $\beta^{k+1}\in\Delta_{\abst{\calV}}$ such that $\supp(\beta^{k+1}) = \calS_{k+1}$.  % namely $\calS_{k+1}$. % and $k:=k+1$. 
\end{enumerate}
\end{algorithmic}%\vspace{-.2cm}
\end{algorithm}

%Before describing our method, let us introduce some notations. We define $\dom F:= \rvA^{-1}(\inter\calK)$ and denote the set of vertices of $\calX$ as $\calV\ne\emptyset$. % In the following, we will sometimes call vertices as ``atoms''. 
%Note that for the implementation of Algorithm~\ref{algo:AFW_SC}, we do not need to know $\calV$ explicitly --- in particular, $\calX$ can be an H-polytope (i.e., $\calX$ is described as a system of linear inequalities). 

 The away-step Frank-Wolfe method for solving~\eqref{eq:poi} is shown in Algorithm~\ref{algo:AFW_SC}.
Let us define $\dom F:= \rvA^{-1}(\inter\calK)$ and write the atomic representation of $\calX$ as $\calX = \conv(\calV)$, where $\calV\subseteq\bbX$ is a finite set of points that contains $\Ext(\calX)$, namely the vertices of $\calX$. We shall call each element in $\calV$ as an {\em atom}. 
Note that we do {\em not} need to have explicit knowledge of $\calV$ in order to run Algorithm~\ref{algo:AFW_SC}. 
%\subseteq\calV$. %of points in $\bbX$
%denote the set of vertices of $\calX$ as $\calV\ne\emptyset$. 
In addition, note that a similar algorithm was proposed in~\cite[Algorithm~8]{Dvu_23}, but there are two subtle but important differences between this algorithm and Algorithm~\ref{algo:AFW_SC}. Specifically,~\cite[Algorithm~8]{Dvu_23} requires that  $x^0\in \dom F\cap \Ext(\calX)$, and define $D_k= \ipt{\nabla^2 F(x^k)d^k}{d^k}^{1/2}$. While this may work for some applications, it can cause serious issues on~\eqref{eq:Dopt}, where $\dom F\cap \Ext(\calX)=\emptyset$ and $D_k$ may be zero for certain feasible $x^k$ (which causes numerical errors when computing $b_k$ in~\eqref{eq:def_alpha_k}).
%  which can be zero since $F$ is degenerate.   % is defined   %First, we only require the starting point

Given a starting point $x^0\in \calX\cap\dom F = \{x\in\calX:\rvA x \in\inter\calK\}$ and its active set of atoms $\calS_0$, we proceed as follows.  At each iteration $k$, we first compute the away atom $v^k$ that solves $\min_{x\in\calV}\;  \lranglet{\nabla F(x^k)}{x}$. Based on $v^k$, we then compute the FW step $d_\rmF^k:= v^k - x^k$ and the {\em FW gap} $G_k$ as in~\eqref{eq:G_k}. Note that by the convexity of $F$, we have that for any $x^*\in\calX^*$, 
\begin{equation}
G_k = \ipt{\nabla F(x^k)}{x^k-v^k}\ge \ipt{\nabla F(x^k)}{x^k-x^*} \ge F(x^k) - F^*:= \delta_k\ge 0.  \label{eq:def_G_k}
\end{equation}
%In addition, 
If $G_k=0$, then we terminate Algorithm~\ref{algo:AFW_SC} and declare  $x^k$ to be an optimal solution of~\eqref{eq:poi}, since $G_k=0$ if and only if $x^k\in\calX^*$. Otherwise, we proceed to compute the descent direction $d^k$. % is an optimal solution of~\eqref{eq:poi} (namely $x^k\in\calX^*$), 

To that end, let $\calS_k\subseteq\calV$ denote set of active atoms at $x^k$, so that $x^k\in\ri\conv(\calS_k)$. If $\abst{\calS_k}=1$, then $x^k\in\calV$ and we simply take $d^k = d^k_\rmF$. Otherwise, we compute the away atom $a^k \in \argmax_{x\in\calS_k} \lranglet{\nabla F(x^k)}{x}$ (by searching over $\calS_k$) and then the away step $d^k_\rmA := x^k - a^k$. We decide between the two directions $d_\rmF^k$ and $d_\rmA^k$ in Step~\ref{item:choose_direction}, by letting $d^k$ be the one that has a higher value evaluated at % the linear function 
$\ipt{-\nabla F(x^k)}{\cdot}$. Note that since $x^k\not\in\calX^*$, we must have 
\begin{equation}
r_k:= \ipt{-\nabla F(x^k)}{d^k}\ge \ipt{-\nabla F(x^k)}{d_\rmF^k} = G_k\ge \delta_k>0, \label{eq:r_k_ge_G_k}
\end{equation}
and hence $d^k$ is indeed a descent direction. After this, we compute the step-size $\alpha_k>0$, by either doing exact line-search or using the adaptive step-size as in~\eqref{eq:def_alpha_k}. Note that for some important applications, including D-optimal design and MHP, the exact line-search procedure can be implemented efficiently --- for details, see~\cite{Khachiyan_96} and~\cite[Appendix~F]{Zhao_23}. 
 After obtaining $d^k$ and $\alpha_k$, we simply set the next iterate $x^{k+1}:= x^k + \alpha_k d^k$.

%An important final step is to 
Finally, we update $\calS_{k+1}$ and $\beta^{k+1}$. % namely the set of active atoms at $x^{k+1}$. 
If $d^k = d^k_\rmA$, then $\beta_{a^k}^{k+1} := \beta_{a^k}^{k} -(1- \beta_{a^k}^{k})\alpha_k  $ and 
$\beta_a^{k+1} := (1+\alpha_k) \beta_a^{k+1}$ for all $a\in\calV\setminus \{a^k\}$. Thus if $\alpha_k = \baralpha_k= {\beta_{a^k}}/({1-\beta_{a^k}})$, then $\beta_{a^k}^{k+1} = 0$ and $\calS_{k+1}:= \calS_{k}\setminus \{a^k\}$; otherwise, we have $\calS_{k+1}:=\calS_{k}$. If $d^k = d^k_\rmF$, then $\beta_{v^k}^{k+1} := (1- \alpha_k)\beta_{v^k}^{k} + \alpha_k$, and $\beta_a^{k+1} := (1-\alpha_k) \beta_a^{k+1}$ for all $a\in\calV\setminus \{v^k\}$. If $v^k\in\calS_k$, then $\calS_{k+1}:=\calS_{k}$; otherwise, we have $\calS_{k+1}:= \calS_{k}\cup \{v^k\}$. Therefore, we have $\abst{\calS_k} -1 \le \abst{\calS_{k+1}}\le \abst{\calS_k} +1$, %at each iteration $k$, 
namely, the size of the active atoms decreases or increases by at most one. For convenience, define $p:=\dim \bbX$. 
When $\abst{\calS_{k+1}}>p+1$ (which typically happens for low to medium dimension $p$), from  Carath\'eodory's theorem and its constructive proof (cf.~\cite[Theorem~17.1]{Rock_70}),  we can remove the redundant atoms in $\calS_{k+1}$ and obtain a smaller set of active atoms $\tilde\calS_{k+1}\subsetneqq \calS_{k+1}$ at $x^{k+1}$ such that $\abst{\tilde\calS_{k+1}}\le p + 1$; for an efficient implementation, see~\cite[Appendix~A]{Beck_17}.  Such a removal of redundant atoms is especially useful  when the size of $\calV$  
%$\abst{\calV}$ 
has a super-linear dependence on $p$ (e.g., $\abst{\calV} = \Theta(p^2)$). 
%an alternative set of active atoms at $x^k$, say $\tilde\calS_{k+1}$, such that $\abst{\tilde\calS_{k+1}}\le \dim\bbX+1$. 

% and  $\calS_k\subseteq \calS_{k+1}$ and $\abst{\calS_k}\ge  \abst{\calS_{k+1}}-1$, namely at most one vertex in $\calS_k$ is ``dropped''. (This happens when $\alpha_k = {\beta_{a^k}}/({1-\beta_{a^k}})$ and hence $\beta^{k+1}_{a^k}=0$). 
%, i.e., $\calS_{k+1}\subseteq \calV$.

\subsection{Some Remarks About Algorithm~\ref{algo:AFW_SC} }

Let us make several remarks regarding the implementation of Algorithm~\ref{algo:AFW_SC}. Recall that $p=\dim \bbX$. 
%For convenience, let us define $p:=\dim \bbX$. 

First, we do not need to know $\calV$ explicitly to implement Algorithm~\ref{algo:AFW_SC}, so long as  %in particular, $\calX$ can be certain H-polytope (namely, it is described by a system of linear inequalities) such that 
we can easily find a feasible starting point $x^0\in\calX\cap\dom F$ and its active set of atoms $\calS_0$. (Some examples will be provided below.)  Specifically, in Step~\ref{item:LMO_a},  if $\calV$ is not explicitly known, we can apply the (primal or dual) simplex method to solve the LP $\min_{x\in\calX}\;  \lranglet{\nabla F(x^k)}{x}$ and obtain an optimal solution $v^k\in\calV$. % and the optimal solution will be a ver
In fact, this approach is preferred even when $\calV$ is explicitly known but has a large size (e.g., $\calX= [0,1]^p$ is a hypercube and $\abst{\calV} = 2^{p}$).

Let us now provide two examples where $\calV$ is not explicitly known, but we can easily identify $x^0$ and $\calS_0$. The first example is the budget-constrained D-optimal design problem~\eqref{eq:Dopt2}, where we have a single budget constraint $c^\top x \le C$. Without loss of generality, let $0\le c_1\le \ldots\le c_m$. If there exists $n\le \barn\le m$ %$C\ge c_{\barn}$ 
such that $c_{\barn}\le C$ and %for some  and there exists $\calI\subseteq [\barn]$ such that 
$\mathsf{span}\{a_i\}_{i\in[\barn]} =\bbR^n$, then  we know that $\{e_i\}_{i\in[\barn]}$ are vertices of the constraint set $\calX=\{x\in\Delta_n:c^\top x\le C\}$, and hence we can take $x^0 = (1/\barn)\sum_{i\in[\barn]} e_i$ and $\calS_0 = \{e_i\}_{i\in[\barn]}$. The second example is the (reformulated) Poisson image de-blurring problem with TV regularization~\eqref{eq:box_TV2}, where we simply observe that $(x^0,r^0) = (Me,Me)$ is a vertex of the constraint set $\calX$ and $\rvA (x^0,r^0)= Ax_0\in\bbR_{++}^N$, and hence we can take $\calS_0 = \{(x^0,r^0)\}$.

Second, note that computing $a^k$ in Step~\ref{item:away_atom_a} is effectively a maximum inner-product search (MIPS) problem, and many efficient (deterministic or randomized) methods have been proposed to solve this problem in recent years (see~\cite{Yu_17MIPS} and references  therein). Of course, when $\calX$ is a unit simplex and $\calV=\{e_i\}_{i=1}^p$ (as in D-optimal design and MHP problems), computing  $a^k$ becomes a trivial task. % can be trivially implemented in $O(p)$ time. 

Third, choosing the adaptive step-size in Step~\ref{item:step_size_SC_a}\ref{item:adapt_step} requires computing the quantity $D_k$. For all the applications described in Sections~\ref{sec:D-opt} and~\ref{sec:two_other}, computing $D_k$  takes $O(p)$ elementary operations, instead of $O(p^2)$. For details, we refer readers to~\cite[Section~2]{Zhao_23}. 

Finally, note that although $\beta^k\in\bbR^{\abst{\calV}}$ (for $k\ge 0$) is stated in %in the description of 
Algorithm~\ref{algo:AFW_SC}, for implementation purpose, we only need to store the nonzero elements of $\beta^k$ and their indices $\calS_k$. This is especially efficient when $\abst{\calS_k}\ll \abst{\calV}$ (e.g., $\calX$ is a hypercube and $\abst{\calV} \ge 2^{p}$). 
%$\abst{\calS_k}$ is much smaller compared to $ \abst{\calV}$.  %, since in certain 
% (namely $l=1$). 

% depending on the structure of $\calX$, we can either perform an exhaustive search over $\calV$, or apply the (primal or dual) simplex method to solve the LP $\min_{x\in\calX}\;  \lranglet{\nabla F(x^k)}{x}$. %, which is typically more efficient 
% Typically, the former approach is preferred if $\calV$ is explicitly known and has a relatively small size (e.g., when $\calX$ is a simplex and $\abst{\calV} = \dim \bbX + 1$), and the latter approach is preferred if $\calV$ is not explicitly known (e.g., $\calX$ is an H-polytope) or has a large size (e.g., $\calX$ is a hypercube and $\abst{\calV} = 2^{\dim \bbX}$). 
%%\vspace{1ex}\noindent 
%{\em 1. Starting point.} 

%Note that the latter approach involves 
%choosing $d^k$, 

\section{Global Linear Convergence of Algorithm~\ref{algo:AFW_SC} } \label{sec:global_lin_conv}

\subsection{Linear Transformation and Characterization of $\calX^*$} \label{sec:char_calX*}

Let us first characterize the optimal solutions $\calX^*$ of~\eqref{eq:poi}. 
To that end, we can rewrite~\eqref{eq:poi} as:
\begin{equation}
(P): \ \ \  {\min}_{x\in\calX} \;[F(x):= \barF(\rvE x)], %\where \rvE x: = (\rvA x, \lranglet{c}{x}) \mbox{ and } \barF(y,t):= f(y) + t
%\barF(\rvA x,\lranglet{c}{x})], 
\end{equation}
where $\rvE x: = (\rvA x, \lranglet{c}{x})\in\bbY\times\bbR$ for all $x\in\bbX$ and
\begin{equation}
\barF(y,t):= f(y) + t, \quad \forall\, y\in\bbY, \;\; \forall t\in\bbR. 
\end{equation}
%for any $y\in\bbY$ and $t\in\bbR$. Indeed, 
%using the variable transformation $z:= (y,t) := (\rvA x,\lranglet{c}{x})$, 
We see that $x^*\in\calX^*$ if and only if %$z^*:= (y^*,t^*) := 
$\rvE x^* =(\rvA x^*,\lranglet{c}{x^*})$ %:= \rvE x^*$ %(\rvA x^*,\lranglet{c}{x^*})$ 
is an optimal solution of %the %can find 
\begin{equation}
{\min} \;\barF(z)\quad \st \quad z\in\rvE(\calX),\quad {\rm where}\;\; z=(y,t)\in \bbY\times\bbR.  \label{eq:poi_yt}
\end{equation}
%where $\calZ:= \{(\rvA x,\ipt{c}{x}):x\in\calX\}$. 
For convenience, let $\calZ:= \rvE(\calX)$. % and $\bbZ:= \bbY\times\bbR$. % in the following.
Since $f$ is strictly convex on its domain, it is easy to show that (see e.g.,~\cite[Lemma~3.1]{Luo_92}) the problem in~\eqref{eq:poi_yt} has a unique solution $z^* := (y^*,t^*)$ and hence 
\begin{equation}
\calX^* =\{x\in\calX:\rvE x = z^*\}=\{x\in\calX:\rvA x = y^*,\lranglet{c}{x} = t^*\}. \label{eq:calX*}
\end{equation}
Of course, the optimality of $z^*$ implies that 
\begin{equation}
\ipt{\nabla f(y^*)}{y-y^*} + (t-t^*) = \ipt{\nabla \barF(z^*)}{z - z^*}\ge 0, \quad\forall\,z\in\calZ. %\rvE(\calX).   
\label{eq:ip_ge_0}
\end{equation}
Note that if $\calY:=\rvA(\calX)$ is a singleton, then~\eqref{eq:poi} reduces to the LP $\min_{x\in\calX} \, c^\top x$.  %namely $\rvA(\calX)= \{y^*\}$, then we can easily find $\calX^*$ from its characterization of in~\eqref{eq:calX*}. 
Hence, in our subsequent discussion, we will always assume that  $\abst{\calY}>1$, and consequently, $\abst{\calZ}\ge\abst{\calY} >1$. %For convenience, %, let us denote $\calZ:= \rvE(\calX)$ in the following.  %Since $\rvA(\calX)$ can be regarded as a ce, we also have  that $\abst{\rvE(\calX)}>1$. 

%For convenience, let us define the linear operator $\rvE: x\mapsto (\rvA x, \lranglet{c}{x})$ and denote $z:= (y,t)$ and $z^* := (y^*,t^*)$. Under this notation, we see that $F(x) = \barF(\rvE x)$ and hence $\nabla F(x) = \rvE^* \nabla \barF(\rvE x)$. 

\subsection{Quadratic Growth Property of $\barF$ on $\calZ$}

Next, let us establish a quadratic growth property of $\barF$ on $\calZ$. To that end, let us define  the radius of $\calY$ with respect to the center $y^*$  and the norm $\normt{\cdot}_{y^*}$ as 
\begin{equation}
R(y^*):= {\sup}_{y\in \calY} \; \normt{y-y^*}_{y^*} \in(0, + \infty). %, \quad \where \calY:=\rvA(\calX), 
\label{eq:def_R_y*}
\end{equation}
and let $\mu_{R(y^*)}>0$ be defined as in~\eqref{eq:quad_omega}. 
%For convenience, let us denote $\calZ:= \rvE(\calX)$, and define %. 
In addition,  define 
\begin{equation}
b_{z^*}:= \max\big\{1,{\max}_{z\in\calZ}\;\ipt{\nabla \barF (z^*)}{z-z^*}\big\}\in(0,+\infty). \label{eq:def_b}
\end{equation}
Based on $\mu_{R(y^*)}$ and $b_{z^*}$, we can define the following symmetric bilinear form $\ipt{\cdot}{\cdot}_{z^*}$ on $\bbY\times\bbR$: for all $z= (y,t)$ and $z'=(y',t')\in \bbY\times\bbR$, we have 
\begin{align}
\ipt{z}{z'}_{z^*}& := \ipt{\nabla^2 \barF(z^*) z}{z'} + \frac{\ipt{\nabla \barF(z^*)}{z}\ipt{\nabla \barF(z^*)}{z'}}{\mu_{R(y^*)} b_{z^*}}\\
&= \ipt{\nabla^2 f(y^*) y}{y'} + \frac{(\ipt{\nabla f(y^*)}{y}+t)(\ipt{\nabla f(y^*)}{y'}+t')}{\mu_{R(y^*)} b_{z^*}}. %\quad \forall\, z= (y,t), z'=(y',t')\in \bbY\times\bbR. 
%\ipt{M_{z^*} z}{z'}, %\quad {\rm where}\quad M_{z^*}:= \nabla^2 \barF(z^*) +  \frac{\nabla \barF(z^*)\nabla \barF(z^*)^\top}{\mu_{R(y^*)} b_{z^*}} . 
\end{align}
Note that  $\ipt{\cdot}{\cdot}_{z^*}$ is positive definite, i.e., $\ipt{z}{z}_{z^*}=0$ if $z=0$ and
$\ipt{z}{z}_{z^*}>0$ if $z\ne 0$. (To see this, consider two cases: i) $y\ne 0$ and ii) $y=0$ but $t\ne 0$.) Therefore, $\ipt{\cdot}{\cdot}_{z^*}$ is an {\em inner product} on $\bbY\times\bbR$, and we let $\normtt{\cdot}_{z^*}$ denote its induced norm, which is defined as
%a nonnegative function $\normtt{\cdot}_{z^*}:\bbY\times\bbR\to\bbR_+$ as follows: 
%distance-like quantity $\rvD(\cdot,{z^*}):\bbY\times \bbR\to\bbR_+$ such that %for any $$, 
\begin{align}
\normtt{z}_{z^*} %&:= \bigg(\normt{y}_{y^*}^2 + \frac{\ipt{\nabla \barF (z^*)}{z}^2}{\mu_{R(y^*)} b_{z^*}}\bigg)^{1/2}\\
& := \bigg(\normt{y}_{y^*}^2 + \frac{({\ipt{\nabla f (y^*)}{y} + t})^2}{\mu_{R(y^*)} b_{z^*}}\bigg)^{1/2}, \quad \forall\, z=(y,t)\in\bbY\times\bbR.  \label{eq:def_normtt_z*}
\end{align}
(Note that we use $\normtt{\cdot}_{z^*}$ to distinguish it from $\normt{\cdot}_{y^*}$, which has a different definition.)
  %we have i) $\normtt{z}_{z^*} = 0$ if and only if $z=0$ and ii) $\normtt{\alpha z}_{z^*} = \abst{\alpha}\normtt{ z}_{z^*}$ for all $z\in\bbY\times \bbR$ and $\alpha\in\bbR$. However, unlike a norm, $\normtt{\cdot}_{z^*}$ does not necessarily satisfy the triangle inequality. For this reason, we shall call $\normtt{\cdot}_{z^*}$ a {\em quasi-norm}. %This 

%\begin{align}
%\hspace{-1ex} \rvD(z,{z^*}) &:= \big(\normt{z-z^*}_{z^*}^2+\mu_{R(y^*)}^{-1} \abst{\lranglet{\nabla \barF(z^*)}{z-z^*}}\big)^{1/2} \label{eq:def_D_z_z^*}\\
%&= \big(\normt{y-y^*}_{y^*}^2+\mu_{R(y^*)}^{-1} \abst{\lranglet{\nabla f(y^*)}{y-y^*} + (t-t^*)}%\big)
%\big)^{1/2}, \;\; \forall\,z=(y,t)\in\bbY\times \bbR,   \label{eq:def_D_z_z^*2}
%\end{align}
%where $\mu_{R(y^*)} > 0 $ is defined in~\eqref{eq:quad_omega}. 
%%Note that $D(z,{z^*})\ge 0$ for all $z\in\rvE(\calX)$ and more importantly, 
%Note that $\rvD(z,{z^*})=0$ if and only if $z=z^*$. 
%$\sigma_{y^*}(\cdot)$ shares some of the properties of a norm, namely, $\sigma_{y^*}(z)\ge 0$ and  $\sigma_{y^*}(z) = \sigma_{y^*}(-z)$  for all $z\in\bbR^{m+1}$ and $\sigma_{y^*}(z) =0$ if and only if $z=0$. 
%The usefulness of $D(\cdot,{z^*})$ will be seen in the following lemma, which establishes an error bound of $\barF$ around $z^*$ on $ \rvE(\calX)$. 

\begin{lemma}[Quadratic growth of $\barF$]\label{lem:eb_F}
%Recall that $z^* = (y^*,t^*)$. 
%For any $z\in\rvE(\calX)$, 
We have 
\begin{equation}
\barF(z)\ge \barF(z^*) + \mu_{R(y^*)}\normtt{z-{z^*}}_{z^*}^2, \quad \forall\,z\in\calZ. 
\end{equation}
\end{lemma}
\begin{proof}%\renewcommand{\qedsymbol}{}
Recall that $z= (y,t)$ and $z^* = (y^*,t^*)$. 
We have
\begin{align}
\barF(z) = f(y) + t&\ge f(y^*) + \ipt{\nabla f(y^*)}{y-y^*} + \omega(\normt{y-y^*}_{y^*}) + t^* + (t-t^*) \label{eq:barF_lb1}\\
&\ge \barF(z^*)  + \mu_{R(y^*)}\normt{y-y^*}_{y^*}^2+ { \ipt{\nabla \barF(z^*)}{z-z^*}}\label{eq:barF_lb2}\\
&\ge \barF(z^*)  + \mu_{R(y^*)}\normt{y-y^*}_{y^*}^2+  b_{z^*}^{-1}\ipt{\nabla \barF(z^*)}{z-z^*}^2\label{eq:barF_lb2.5}\\
&= \barF(z^*) + \mu_{R(y^*)}\normtt{z-{z^*}}_{z^*}^2,\label{eq:barF_lb3}
\end{align}
where in~\eqref{eq:barF_lb1} we use Lemma~\ref{lem:lb_f}, in~\eqref{eq:barF_lb2} we use $\normt{y-y^*}_{y^*}\le R(y^*)$ for any $y\in\rvA(\calX)$ and Lemma~\ref{lem:omega}, in~\eqref{eq:barF_lb2.5} we use the optimality condition in~\eqref{eq:ip_ge_0} and $b_{z^*}\ge \ipt{\nabla \barF(y^*)}{z-z^*} $ for all $z\in\calZ$, % as well as ~\eqref{eq:ip_ge_0},  
and finally, in~\eqref{eq:barF_lb3} we use the definition of $\normtt{\cdot}_{z^*}$ in~\eqref{eq:def_normtt_z*}. 
%\begin{equation}
%\ipt{\nabla f(y^*)}{y-y^*} + (t-t^*) = \ipt{\nabla \barF(z^*)}{z - z^*}\ge 0, \quad\forall\,z\in\rvE(\calX).\nn\tag*{\qedhere} %\qedhere  \hspace{3ex}\qed
%\end{equation}
\end{proof}
\subsection{A Condition Number About $\calX$ and $\calX^*$}

As per discussion in Section~\ref{sec:char_calX*}, we focus on the case where $\abst{\calZ}>1$, and hence $\calX\setminus\calX^*\neq\emptyset$. Consequently, fix any $x\in\calX\setminus\calX^*$ and define $z:= \rvE x\ne z^*$. %From~\eqref{eq:calX*}, we know that $z\neq z^*$. 
Accordingly, let us define %$$d:= \frac{z - z^*}{\rvD( z, z^*)}.$$
\begin{align}
d:= \frac{z - z^*}{\normtt{ z- z^*}_{z^*}} \ne 0 \quad \mbox{and}\quad \Phi_\calS(x,\calX^*) := \min_{p:\ipt{p}{d}=1}\; \max_{a\in\calS,v\in\calV}\; \ipt{p}{\rvE a- \rvE v}, \label{eq:def_Phi_S_x_x*}
\end{align}
where $\calS\subseteq \calV$ satisfies  that $x\in\ri \conv(\calS)$, and we call $\calS$ an atomic representation of $x$.  Let ${\sf AR}(x)$ denote the set of all atomic representations at $x$, namely ${\sf AR}(x):=\{\calS\subseteq\calV:x\in\ri\conv(\calS)\}$. Then, we can define the following important geometric constant: 
\begin{equation}
\Phi(\calX,\calX^*):= {\min}_{x\in\calX\setminus\calX^*}\;\;{\min}_{\calS\in {\sf AR}(x)}\;\;\Phi_\calS(x,\calX^*). \label{eq:def_Phi_X_X*}
\end{equation}

%${\min}_{\calS\in\calS(x)}\;$
%In addition, Based on $d$ and ${\sf AR}(x)$, let us define 
% % and any $x^*\in\calX^*$, let us define the following quantity: % as in~\cite{Pena_19}:
%\begin{equation}
%\Phi(x,x^*) := {\min}_{\calS\in\calS(x)}\;{\min}_{p:\ipt{p}{d}=1}\; {\max}_{a\in\calS,v\in\calV}\; \ipt{p}{\rvE a- \rvE v}, \label{eq:def_Phi_x_x*} %\quad \forall\,x\in\calX,\,\forall\,x^*\in\calX^*,
%\end{equation}
%where $\scS(x):= \{\calS\subseteq\calV:x\in\ri\conv(\calS)\}$, namely the collection of all the set of ``active'' vertices at $x\in\calX\setminus\calX^*$ %$\calS$ such that denotes a set of ``active'' vertices at $x$ such that 
%and $$d:= \frac{\rvE x- \rvE x^*}{D(\rvE x,\rvE x^*)}.$$ (Note that since $x\in\calX\setminus\calX^*$, $\rvE x\ne  \rvE x^*$ and $D(\rvE x,\rvE x^*)>0$.) We then define 

%Using the same argument as in the proof of~\cite[Proposition~3]{Pena_19} (by replacing the Euclidean norm $\normt{\cdot}$ with the local norm $\normt{\cdot}_{y^*}$), 

Next, let us show that the constant $\Phi(\calX,\calX^*)$ can be lower bounded by certain ``facial distance'' of the polytope $\calZ$, which is always positive. 

%%The following lemma provides a lower bound on $\Phi(\calX,\calX^*)$ in terms of the ``facial distance'' of the polytope $\rvE(\calX)$. The statement of this lemma requires defining the following quasi-norm $\normt{\cdot}_{\rvq,z^*}$ (which satisfies all the properties of a norm except triangle inequality):
%%\begin{align*}
%%\normt{z}_{\rvq,z^*} &:= \big(\normt{y}_{y^*}^2 + \abst{\ipt{\nabla \barF (z^*)}{z}}^2/(\mu(R(y^*)) g)\big)^{1/2}, \\
%%& = \big(\normt{y}_{y^*}^2 + \abst{\ipt{\nabla f (y^*)}{y} + t}^2/(\mu(R(y^*)) g)\big)^{1/2}, \qquad \forall\, z=(y,t)\in\bbR^{m+1},
%%\end{align*}
%%where
%%\begin{equation}
%%g:= \max\{1,\barg\}\ge 1 \andd \barg:= {\max}_{z\in\rvE(\calX)}\;\ipt{\nabla \barF (z^*)}{z-z^*}. \label{eq:def_g}
%%\end{equation}

\begin{lemma} \label{lem:lb_Phi}
%If $\rvE(\calX)\ne\{z^*\}$ (namely $\rvE(\calX)$ is not a singleton), then 
Let $\calF_{z^*}$ denote the minimal face of $\calZ $ that contains $z^*$, and let ${\sf Faces}(\calF_{z^*})$ denote the set of faces of $\calF_{z^*}$. Then we have 
\begin{align*}
\Phi(\calX,\calX^*) \ge {\min}_{\emptyset\ne \calF\in {\sf Faces} (\calF_{z^*}), \,\calF\ne \calZ}\; \dist_{\normtt{\cdot}_{z^*}}\big(\conv(\rvE(\calV)\setminus \calF), \calF\big) > 0, 
\end{align*}
where $\dist_{\normtt{\cdot}_{z^*}}(\calA, \calB):= \inf_{z\in\calA,z'\in\calB}\; \normtt{z-z'}_{z^*}$ for nonempty sets $\calA,\calB\subseteq\bbY\times \bbR$.  
%where $\calF(z^*)$ denote the minimal face of $\rvE(\calX)$ that contains $z^*$, $\scF(\calF(z^*))$ denotes the set of faces of $\calF(z^*)$ and % and \\
%$$\dist_{\normt{\cdot}_{\rvq,z^*}}(\calG,\conv(\rvE(\calV)\setminus\calG)):= \min_{u\in\calG,\,x\in \conv(\rvE(\calV)\setminus\calG)}\,\normt{x-u}_{\rvq,z^*}.$$
\end{lemma}

\begin{proof}
See Appendix~\ref{app:positivity_Phi}. 
\end{proof}

%In Appendix~\ref{app:positivity_Phi}, using similar arguments as in~\cite{Pena_19},  we show that $\Phi(\calX,\calX^*)>0$ so long as $\rvE(\calX)\ne\{z^*\}$, namely $\rvE(\calX)$ is not a singleton. 

The importance of the condition number $\Phi(\calX,\calX^*)$ %was named ``local facial distance'' in~\cite{Pena_19}, and 
%plays an important role
is shown in the following lemma.

\begin{lemma}\label{lem:lb_innerp}
For any optimal solution $x^*\in\calX^*$, any %sub-optimal solution 
$x\in\calX\setminus\calX^*$ and any $\calS\in{\sf AR}(x)$, %atomic representation $\calS$ of $x$ (i.e., $\calS\in{\sf AR}(x)$), 
we have 
\begin{equation}
{\max}_{a\in\calS,v\in\calV}\;\ipt{\nabla F(x)}{a- v}\ge \Phi(\calX,\calX^*) \frac{\ipt{\nabla F(x)}{ x-  x^*}}{\normtt{\rvE x - z^*}_{z^*}}. \label{eq:lb_Phi_X_X*}
\end{equation}
\end{lemma}

\begin{proof}
Since $x^*\in\calX^*$, we know $z:= \rvE x\ne z^*$ and $\ipt{\nabla F(x)}{ x-  x^*}\ge F(x) - F(x^*)>0$. In addition, let $d\ne 0$ be defined as in~\eqref{eq:def_Phi_S_x_x*}, and we have
\begin{equation}
\ipt{\nabla \barF(z)}{d} = \frac{\ipt{\nabla \barF(z)}{z-z^*}}{\normtt{z - z^*}_{z^*}}=\frac{\ipt{\nabla F(x)}{x-x^*}}{\normtt{\rvE x - z^*}_{z^*}}>0. 
\end{equation}
%As a result, 
Since $\ipt{\nabla F(x)}{a- v}= \ipt{\nabla \barF(z)}{\rvE a- \rvE v}$, we can write~\eqref{eq:lb_Phi_X_X*}  equivalently as
\begin{equation}
{\max}_{a\in\calS,v\in\calV}\;\ipt{\barp}{\rvE a- \rvE v}\ge \Phi(\calX,\calX^*), \quad \mbox{where}\quad \barp:= {\nabla \barF(z)}/{\ipt{\nabla \barF(z)}{d}}. \label{eq:lb_Phi_X_X*_equiv}
\end{equation}
Since $\ipt{\barp}{d}=1$, from the definitions of $\Phi_\calS(x,\calX^*)$ and $\Phi(\calX,\calX^*)$ (in~\eqref{eq:def_Phi_S_x_x*} and~\eqref{eq:def_Phi_X_X*}, respectively), we know that~\eqref{eq:lb_Phi_X_X*_equiv} clearly holds. %This completes the proof.  %we have ${\max}_{a\in\calS,v\in\calV}\;\ipt{p}{\rvE a- \rvE v}\ge \Phi_\calS(x,\calX^*)$ and 
\end{proof}

%Equipped with Lemma~\ref{lem:lb_innerp}, we are ready to present the global linear convergence result of Algorithm~\ref{algo:AFW_SC}. %Before doing so, let us 

\subsection{Global Linear Convergence Rate}

%For convenience, 
Let us first define the objective gap at any feasible point $x\in\calX$ as 
\begin{equation}
\delta(x):= F(x) - F^*,  \label{eq:delta}
\end{equation}
and define the shorthand notation $\delta_k:= \delta(x^k)$, for all $k\ge 0$. In this section, we will state the global linear convergence rate of Algorithm~\ref{algo:AFW_SC} in terms of both i) the sequence of objective gaps $\{\delta_k\}_{k\ge 0}$ and ii) the sequence $\{r_k\}_{k\ge 0}$. Note that 
since $r_k\ge G_k$ for all $k\ge 0$, %we can also obtain
 the global linear convergence rate of $\{r_k\}_{k\ge 0}$ also holds for the sequence of FW gaps $\{G_k\}_{k\ge 0}$. %The rate 
 Before stating the convergence rates, we need to define the following quantities. The first one is 
\begin{equation}
B:={\max}_{x,x'\in\calX}\;\ipt{c}{x-x'}\ge 0, \label{eq:def_B}
\end{equation}
namely the variation of the linear function $\ipt{c}{\cdot}$ on $\calX$. Note that for any norm $\normt{\cdot}$ on $\bbX$ (with dual norm $\normt{\cdot}_*$), we can upper bound $B$ as  %$B\le \normt{c}_* D_{\calX,\normt{\cdot}}$, %where $D_{\calX,\normt{\cdot}}: = {\max}_{x,x'\in\calX}\;\normt{x-x'}$
\begin{align*}
B\le \normt{c}_* \Diam_{\normt{\cdot}}(\calX), \quad \mbox{where} \quad \Diam_{\normt{\cdot}}(\calX): = {\max}_{x,x'\in\calX}\;\normt{x-x'}. 
\end{align*}
%denotes the diameter of $\calX$ under $\normt{\cdot}$. 
The second quantity is
\begin{equation}
q:= \min\{\abst{\calW}:\; \calW\subseteq \calV\;\;\mbox{such that}\;\; \conv{\calW}\cap\dom F\ne \emptyset\}\ge 1,  
\end{equation}
where recall that $\dom F:=\rvA^{-1}(\inter\calK)$. %\{x\in\bbX:\rvA x \in\inter\calK\}$.  %$\dom F= \rvA^{-1}(\inter\calK):= $
In words, $q$ is the minimum number of atoms needed to ``produce'' (via convex combination) a feasible point of~\eqref{eq:poi}. 
%Finally, %for any feasible point $x\in\calX$, 
%we define the objective gap at any feasible point $x\in\calX$ as 
%\begin{equation}
%\delta(x):= F(x) - F^*, 
%\end{equation}
%and define the shorthand notation $\delta_k:= \delta(x^k)$, for all $k\ge 0$. 
%, which we denote by $q\ge 1$, 
%such that there exists a convex combination of these atoms that lies in $\dom f$. 

\begin{theorem}[Global linear convergence of $\{\delta_k\}_{k\ge 0}$] \label{thm:global_lin_conv}
%Let $B:=\max_{x,x'\in\calX}\;\ipt{c}{x-x'}\ge 0$ denote  
In Algorithm~\ref{algo:AFW_SC}, no matter exact line-search or adaptive step-size is used in Step~\ref{item:step_size_SC_a}, the following holds: 
\begin{enumerate}[label = (\roman*),leftmargin = 2em]
\item  \label{item:monotone}
The sequence of objective gaps $\{\delta_k\}_{k\ge 0}$ is strictly decreasing (until termination).
\item \label{item:lin_rate}
For all $k\ge 0$, define $k_{\rm eff}:= \ceil{\max\{(k-\abst{\calS_0}+q)/2,0\}}$, and we have
\begin{equation}
\delta_k\le (1-\rho)^{k_{\rm eff}}\delta_0, \quad\mbox{where} \quad\rho := \min\left\{\frac{1}{5.3(\delta_0 + \theta + B)}, \frac{\mu_{R(y^*)} \Phi(\calX,\calX^*)^2}{42.4 (\theta + B)^2}\right\}.  \label{eq:global_lin_rate}
\end{equation}
\end{enumerate}
\end{theorem}

\noindent 
To  state the convergence rate of $\{r_k\}_{k\ge 0}$, we need to define the following affine-invariant quantity: %\vspace{-1ex}
\begin{equation}
%\hspace{-5ex}
\begin{split}
&\bar D:= \max\big\{\ipt{\nabla^2 F(x) d}{d}^{1/2}:\;\;x\in\calS_F(x^0),\; d\in\calX-\calX\big\}<+\infty,\quad\\   
\where  &\calS_F(x^0):= \{x\in\calX:F(x)\le F(x^0)\} \;\;\andd\;\; \calX-\calX:= \{x-x':x,x'\in\calX\}. 
 % < +\infty,%{\max}_{y\in\rvA(\calS_F(x^0))}\;\;{\max}_{d\in\calY-\calY} \;\; 
\end{split}%\nn%\label{eq:barD}
\label{eq:barD}
\end{equation}
%where $\calS_F(x^0):= \{x\in\calX:F(x)\le F(x^0)\}$, % denotes the sub-level set of $F$ at $x^0$ w.r.t.\ $\calX$, $\calY:= \rvA(\calX)$ %is a polytope and $\calY-\calY:= \{y-y':y,y'\in\calY\}$. 
Clearly, from the definition of $D_k$ in Step~\ref{item:step_size_SC_a}\ref{item:adapt_step} of Algorithm~\ref{algo:AFW_SC},  we have
\begin{equation}
D_k\le \bar D, \quad \forall\,k\ge 0. \label{eq:D_k_le_D}
\end{equation}
To see $\bar D<+\infty$, first note that $\calS_F(x^0)\ne \emptyset$ is compact (since $\calX$ is bounded and $F$ is closed)  and $\rvA(\calS_F(x^0))\subseteq\inter\calK$. Since $\nabla^2 f(\cdot)$ is continuous on $\inter\calK$ and $\calY$ is compact, we have $\bar D<+\infty$. 
%Then, This, together with the continuity of $\nabla^2 f(\cdot)$ on $\inter\calK$ and the compactness of $\calY$, implies that $\bar D<+\infty$.  

%To start with, observe that the sequence $\{D_k\}_{k\ge 0}$ can be uniformly upper bounded by the 

\begin{theorem}[Global linear convergence of $\{r_k\}_{k\ge 0}$] \label{thm:global_lin_conv_r_k}
In Algorithm~\ref{algo:AFW_SC}, no matter exact line-search or adaptive step-size is used in Step~\ref{item:step_size_SC_a}, we have %for all $k\ge 0$, %the following holds: 
\begin{equation}
r_k\le 
\left\{
	\begin{array}{ll}
		4(1-\rho)^{k_{\rm eff}}\delta_0\max\{\bar D, 1\},  & \mbox{if } \delta_k >1 \\[1.5ex]
		4\sqrt{1-\rho}^{k_{\rm eff}}\sqrt{\delta_0}\max\{\bar D, 1\}, & \mbox{if } \delta_k \le 1
	\end{array}
\right.  , \quad \forall\,k\ge 0, \label{eq:conv_r_k}
\end{equation}
where $k_{\rm eff}$ and $\rho$ are defined in Theorem~\ref{thm:global_lin_conv}. 
\end{theorem}

\begin{remark}\label{rmk:1/2-worse}
Essentially, Theorem~\ref{thm:global_lin_conv_r_k} states that $\{r_k\}_{k\ge 0}$ (and hence $\{G_k\}_{k\ge 0}$) converges linearly at a rate that is about (1/2)-order of magnitude worse compared to that of $\{\delta_k\}_{k\ge 0}$. 
\end{remark}

Before proving Theorems~\ref{thm:global_lin_conv} and~\ref{thm:global_lin_conv_r_k}, %(in Sections~\ref{sec:proof_delta_k} and~\ref{sec:proof_r_k}, respectively), 
let us first make a few remarks about them.

\subsubsection{Remarks on Theorems~\ref{thm:global_lin_conv} and~\ref{thm:global_lin_conv_r_k} }

%Before proving Theorem~\ref{thm:global_lin_conv}, let us make a few remarks about its part~\ref{item:lin_rate}. 

First, $k_{\rm eff}$ denotes the number of ``effective'' iterations within the first $k$ iterations --- we call an iteration ``effective'' if within this iteration, the objective gap is reduced by a factor of at least $1-\rho$. Clearly, we have $k_{\rm eff}\ge 1$ if and only if $k\ge \abst{\calS_0}-q+1$.

Second, for all of the applications in Sections~\ref{sec:D-opt} and~\ref{sec:two_other},  %some important applications including D-optimal design and MHP,
\cal we can obtain a starting point $x^0$ such that $\abst{\calS_0}=q$ fairly easily. For the purpose of discussion, let $\calV=\Ext(\calX)$. 
%Specifically, f
For the D-optimal design problem~\eqref{eq:Dopt}, since $q = n$, it suffices to choose $n$ linear independent points $\calA_\calI:=\{a_i\}_{i\in\calI}$ from $\calA$ (with $\calI\subseteq[n]$ and $\abst{\calI}=n$)  and let $x^0=(1/n)\sum_{i\in\calI} e_i$. Selecting $\calA_\calI$ can be accomplished simply by Gaussian elimination on the $n\times m$ matrix $[a_1\;\cdots\;a_m]$, which takes $O(mn^2)$ elementary operations. For the MHP problem in~\eqref{eq:L1_MLE_rewrite}, we can simply take $(\mu^0,a^0) = (1,0)$, which is a vertex of $\Delta_{m+1}$, and hence $q=1$. Similarly, for the Poisson image de-blurring problem with TV regularization in~\eqref{eq:box_TV2}, we can take $(x^0,r^0) = (Me,Me)$, which a vertex of the constraint set and hence $q=1$. 
 
Third, note that the global linear convergence rate in~\eqref{eq:global_lin_rate} is {\em affine invariant}. %Specifically, 
This is because in~\eqref{eq:global_lin_rate}, $\rho$  only depends on five quantities, namely $\delta_0$, $\theta$, $B$, $R(y^*)$ and $\Phi(\calX,\calX^*)$.  (Note that $\mu_{R(y^*)}$ is simply a univariate function of $R(y^*)$; cf.~\eqref{eq:quad_omega}.) The first four quantities are easily seen to be affine invariant, and upon careful examination, we also see that $\Phi(\calX,\calX^*)$ is affine invariant (cf.~\eqref{eq:def_Phi_X_X*}). In addition, since $\bar D$ is affine-invariant (cf.~\eqref{eq:barD}), the convergence rate of $\{r_k\}_{k\ge 0}$ in~\eqref{eq:conv_r_k} is also affine-invariant. 
%Therefore, 
%we see that 
%the last one (i.e., $\Phi(\calX,\calX^*)$) is affine-invariant. 
 
%Third, if we choose $x^0$ such that $\abst{\calS_0}=q$, then the iteration complexity of Algorithm~\ref{algo:AFW_SC} to obtain a $\varepsilon$-optimal point $x\in\calX$ (namely, $\delta(x)\le \varepsilon$) is 
%\begin{equation}
%O\bigg(\max\bigg\{\delta_0 + \theta + B,\;\; \frac{(\theta + B)^2}{\mu_{R(y^*)} \Phi(\calX,\calX^*)^2}\bigg\}\ln\bigg(\frac{\delta_0}{\varepsilon}\bigg)\bigg). 
%\end{equation}
%an $\varepsilon$ objective gap 
 
%
% involved in the 
 %, where  and $\{a_i\}{i\in\calI}$ are linear independent. 

%, we mean the iterations % steps denote those 
%let us collect a few lemmas from~\cite{Zhao_23} below that will be useful in our proof. 

%Next, let us prove Theorem~\ref{thm:global_lin_conv}. 

\subsubsection{A Useful Lemma}

To prove Theorem~\ref{thm:global_lin_conv}, we need the following useful lemma. % is important for proving Theorem 

\begin{lemma}\label{keyinequality} 
For all $k \ge 0$ it holds that 
\begin{equation}\label{eq:bound_Dk}
D_k \le \max\{r_k,\sqrt{\theta}\} +\theta+ B\ . 
\end{equation}
\end{lemma}

\noindent 
Our proof use the following properties of the function class $\calB_\theta(\calK)$,  in addition to those listed  in Section~\ref{sec:LHSCB}. These properties are listed below. 

\begin{lemma}[{see Nesterov and Nemirovskii~\cite[Corollary~2.3.1, Proposition~2.3.2, Proposition~2.3.4, and Corollary 2.3.3]{Nest_94}}]\label{lem:LHSCB}
If $f\in \calB_\theta(\calK)$, %be a $\theta$-LHSCB defined on a nonempty, closed and convex cone $\calK$. %, where $\theta\ge 1$. 
then for any $y\in\inter\calK$, % and $u\in\bbR^n$,
 we have %$h$ satisfies that 
\begin{enumerate}[start = 5,label = {\rm (P\arabic*)}, leftmargin = 3.9em]
\item  \label{item:def_SCB} $\abs{\lranglet{\nabla f(y)}{u}} \le \sqrt{\theta}\|u\|_y$ $\quad \forall\,u\in\bbY$, %\;\forall\,x\in\inter\calK, \;\forall\,u\in\bbR^n.$
\item \label{item:recession_cone} $\|v\|_y \le -\lranglet{\nabla f(y)}{v}$ $\quad \forall\,v\in\calK$,
\item \label{item:ub_innerp} $\ipt{\nabla f(y)}{v-y}\le \theta$ $\quad \forall\, v\in\calK$,
\item \label{item:Hessian_grad} $\nabla^2 f(y)y = -\nabla f(y)$,
%$\lranglet{\nabla f(y)}{u} = -\lranglet{H(u)u}{w}$ $\quad \forall\,w\in\bbR^m$,
\item \label{item:grad_identity} $\lranglet{\nabla f(y)}{y} = -\theta$. %and
%\item \label{item:thetageone} $\theta \ge 1$. \qed
\end{enumerate}
\end{lemma}
\noindent Note that as a direct consequence of~\ref{item:Hessian_grad} and~\ref{item:grad_identity}, we have
\begin{equation}
\normt{y}_y^2 = \theta,\quad \forall\,y\in\inter\calK.  \label{eq:norm_y^2}
\end{equation}

%To show Lemma~\ref{keyinequality},  

\begin{proof}[Proof of Lemma~\ref{keyinequality}]
Let us consider two cases, namely $d^k = d^k_\rmF$ and $d^k = d^k_\rmA$. 
When $d^k = d^k_\rmF$, we have  $r_k = G_k$, and from~\cite[Proposition 2.3]{Zhao_23}, we know that  
\begin{equation}
D_k \le G_k +\theta+ B = r_k +\theta+ B.  \label{eq:case_FW_step}
\end{equation}
%$$. Thus~\eqref{eq:bound_Dk} clearly holds. 
When $d^k = d^k_\rmA$, we have 
\begin{equation}
r_k = \ipt{\nabla f(y^k)}{w^k-y^k} + \ipt{c}{a^k-x^k}\label{eq:def_r_k}
\end{equation}
and 
\begin{align}
D_k^2= \normt{A(x^k-a^k)}^2_{y^k}  = \normt{y^k-w^k}^2_{y^k}&= \normt{y^k}^2_{y^k} - 2\ipt{\nabla^2 f(y^k) y^k}{w^k} + \normt{w^k}^2_{y^k}\nn\\
 &= 2\ipt{\nabla^2 f(y^k) y^k}{y^k-w^k} + \normt{w^k}^2_{y^k} - \theta \label{eq:ki_1}\\
 &= 2\ipt{\nabla f(y^k) }{w^k-y^k} + \normt{w^k}^2_{y^k} - \theta\label{eq:ki_2}\\
 &\le \theta + \ipt{\nabla f(y^k)}{w^k}^2,\label{eq:ki_3}
\end{align}
where in~\eqref{eq:ki_1} we use~\eqref{eq:norm_y^2}, in~\eqref{eq:ki_2} we use~\ref{item:Hessian_grad}, and in~\eqref{eq:ki_3} we use $w^k\in\calK$,~\ref{item:ub_innerp} and~\ref{item:recession_cone}.  Next, let us show 
\begin{align}
-(\theta + B)\le \ipt{\nabla f(y^k)}{w^k}\le 0. \label{eq:lb_ip}
\end{align}
%In fact, the second inequality is 
Indeed, from~\ref{item:grad_identity}, we know that $\ipt{\nabla f(y^k)}{w^k}= \ipt{\nabla f(y^k)}{w^k-y^k}- \theta$. Using~\ref{item:ub_innerp}, we conclude that $\ipt{\nabla f(y^k)}{w^k}\le 0$. In addition, 
 from~\eqref{eq:def_r_k} and~\eqref{eq:def_B}, we have
 \begin{align}
 \ipt{\nabla f(y^k)}{w^k}  = r_k - \ipt{c}{a^k-x^k} - \theta\ge -(B+\theta). 
 \end{align}
 Now, from~\eqref{eq:ki_3} and~\eqref{eq:lb_ip}, we know that $D_k^2\le \theta + (\theta+B)^2$ and hence
\begin{equation}
D_k\le \sqrt{\theta} + \theta+B. \label{eq:case_away_step}
\end{equation} 
   %\le 2\theta + B$ (since $\theta\ge 1$). 
By combining~\eqref{eq:case_FW_step} and~\eqref{eq:case_away_step}, we complete the proof. 
  %in this case.  %\le r_k +2\theta+ B$. %Since $r_k = G_k$ in this case, we k
\end{proof}

%\begin{proof}
%See Appendix~\ref{app:proof_keyineq}. 
%\end{proof}

\subsubsection{Proof of Theorem~\ref{thm:global_lin_conv}} \label{sec:proof_delta_k}

%\begin{proof}[Proof of Theorem~\ref{thm:global_lin_conv}]
Note that it suffices to prove Theorem~\ref{thm:global_lin_conv} for Algorithm~\ref{algo:AFW_SC} with adaptive step-size (cf.\ Step~\ref{item:step_size_SC_a}\ref{item:adapt_step}), since at each iteration, exact line-search (cf.\ Step~\ref{item:step_size_SC_a}\ref{item:exact_LS}) will no less progress in reducing the objective gap compared to adaptive step-size. To start, note that at any iteration $k\ge 0$ such that $G_k>0$,  we have $x^k\in\calX\setminus\calX^*$. Since $r_k\ge G_k$, we also have $r_k>0$.  
Let us first note that 
\begin{equation}
0\le \alpha_k D_k<1, \quad \forall\,k\ge 0. \label{eq:alpha_k_D_k}
\end{equation}
%$$ for all $k\ge 0$. 
Indeed, this is obvious if $D_k=0$; otherwise, from~\eqref{eq:def_alpha_k}, we have $\alpha_k D_k\le b_k D_k = {r_k}/({r_k+D_k}) < 1. $ %As a result, we 
Since $D_k = \normt{\rvA d^k}_{y^k}$, from~\eqref{eq:alpha_k_D_k}, we have $y^{k}+\alpha_k \rvA d^k\in\calB_{y^k}(y^k,1)$. Consequently, from Lemma~\ref{lem:ub_f}, we have 
%\begin{equation}
%\alpha_k D_k\le b_k D_k = \frac{r_k}{r_k+D_k} < 1. 
%\end{equation}Note that we have
\begin{align*}
F(x^{k+1}) = F(x^{k}+\alpha_k d^k) &= f(y^{k}+\alpha_k \rvA d^k)+\lranglet{c}{x^{k}+\alpha_k d^k}\\
&\le f(y^k) + \alpha_k \ipt{\nabla f(y^k)}{\rvA d^k} + \omega^*\big(\alpha_k \normt{\rvA d^k}_{y^{k}}\big) +\lranglet{c}{x^{k}}+\alpha_k\ipt{c}{ d^k}\\
&= F(x^k) - \alpha_k \ipt{-\nabla F(x^k)}{d^k}  + \omega^*\big(\alpha_k \normt{\rvA d^k}_{y^{k}}\big)\\
&= F(x^k) - \alpha_k r_k  + \omega^*(\alpha_k D_k).\nt\label{eq:curv_alpha_k}
\end{align*}
To further proceed, we shall consider the full-step and non-full-step cases separately. \\%[1ex]

\noindent 
\underline{The full-step case: $\alpha_k = \baralpha_k$.} In this case, we have two possibilities, namely  $D_k =0$ and $D_k >0$. If $D_k =0$, by~\eqref{eq:curv_alpha_k} and the fact that $r_k\ge G_k\ge \delta_k$, we have
\begin{align}
F(x^{k+1})\le F(x^k) - \baralpha_k r_k\le F(x^k) - \baralpha_k \delta_k\quad\Longleftrightarrow\quad \delta_{k+1}\le (1-\baralpha_k)\delta_k. \label{eq:delta_k_D_k=0}
\end{align} 
If $D_k >0$, then we have $b_k\ge \baralpha_k$, which amounts to 
\begin{equation}
 (1-\baralpha_k D_k)r_k \ge \baralpha_k D_k^2>0.
\end{equation}
%As a result, 
Since $\alpha_k = \baralpha_k$ and $\alpha_k D_k<1$ (cf.~\eqref{eq:alpha_k_D_k}), we have $\baralpha_k D_k<1$ and hence 
\begin{equation}
r_k \ge \frac{\baralpha_k D_k^2}{1-\baralpha_k D_k}. \label{eq:rk_lb}
\end{equation}
%In sub-case (i), using~\eqref{eq:curv_alpha_k} and that $r_k\ge G_k\ge \delta_k$, we have
%\begin{align*}
%F(x^{k+1})\le F(x^k) - \baralpha_k r_k\le F(x^k) - \baralpha_k \delta_k\quad\Longleftrightarrow\quad \delta_{k+1}\le (1-\baralpha_k)\delta_k. 
%\end{align*} 
%In sub-case (ii), we have $\baralpha_k D_k<1$ and As a result, 
In addition, by~\eqref{eq:curv_alpha_k} and Lemma~\ref{lem:ub_omega*}, we have
\begin{align*}
F(x^{k+1}) &\le F(x^k) - \baralpha_k r_k + \omega^*(\baralpha_k D_k)\le F(x^k) - \baralpha_k r_k + \frac{\baralpha_k^2 D_k^2}{2(1-\baralpha_k D_k)}. 
\end{align*}
This, together with~\eqref{eq:rk_lb} and the fact that $r_k\ge G_k\ge \delta_k$, implies that % and~\eqref{eq:rk_Gk_delta_k},   we have
\begin{align}
F(x^{k+1})\le F(x^k) - \frac{\baralpha_k r_k}{2}\le F(x^k) - \frac{\baralpha_k \delta_k}{2} \quad\Longleftrightarrow\quad \delta_{k+1} \le (1-\baralpha_k/2) \delta_k. \label{eq:delta_k_D_k>0}
\end{align}
%where in the last step we use $r_k\ge G_k\ge \delta_k$. %This is equivalent to 
Thus, from~\eqref{eq:delta_k_D_k=0} and~\eqref{eq:delta_k_D_k>0}, we know that no matter $D_k=0$ or $D_k>0$, we always have 
\begin{equation}
\delta_{k+1} \le (1-\baralpha_k/2) \delta_k, \quad \forall\,k\ge 0. 
\end{equation}
If $\baralpha_k=1$, then $\delta_{k+1}\le \delta_{k}/2$; otherwise, we have $\baralpha_k = \frac{\beta_{a^k}}{1-\beta_{a^k}}$ and $\delta_{k+1}< \delta_{k}$. Note that in the latter case %Algorithm~\ref{algo:AFW_SC} may not 
we may not have a sufficient decrease of the objective gap;  however, since this corresponds to a ``drop'' step, it happens at most $\floor{(\abst{\calS_0}+k-q)/2}$ times within the first $k$ steps. % (see Section~\ref{sec:sublinear}). 
(To see this, let $k_\rmd$ denote the number of  ``drop'' steps within the first $k$ steps, and note that  we always have $\abst{\calS_0}+(k-k_\rmd) - k_\rmd\ge q$.)\\%[1ex]

\noindent 
\underline{The non-full-step case: $\alpha_k < \baralpha_k$.} In this case, we have $D_k>0$ and $\alpha_k=b_k\le \baralpha_k$.  By substituting $\alpha_k = b_k$ into~\eqref{eq:curv_alpha_k}, and using Lemma~\ref{keyinequality} and the monotonicity of $\omega^*$ on $\bbR_+$, we have  
%Prop.\ 2.3 in~\cite{Zhao_23}, we have $D_k \le G_k + \theta\le r_k + \theta$, and so 
\begin{equation}
\delta_{k+1}\le \delta_k - \omega\left(\frac{r_k}{D_k}\right) \le  \delta_k - \omega\left(\frac{r_k}{\max\{r_k,\sqrt{\theta}\} + \theta + B}\right). \label{eq:decsent2}
\end{equation}
Since $r_k>0$, we have $\delta_{k+1}< \delta_k$. This, together with the discussion in the full-step case, shows part~\ref{item:monotone} of Theorem~\ref{thm:global_lin_conv}. Next, let us consider two cases, namely  $r_k> \theta + B$ and  $r_k\le \theta + B$. 
If $r_k> \theta + B$, then $r_k\ge \sqrt{\theta}$ and %by the 
\begin{align*}
\frac{r_k}{\max\{r_k,\sqrt{\theta}\} + \theta + B} = \frac{r_k}{r_k + \theta + B} > \frac{1}{2}, 
\end{align*}
%${r_k}/({r_k + \theta + B})>1/2$ 
and by~\eqref{eq:lin_omega} in Lemma~\ref{lem:omega}, we have 
\begin{equation}
\omega\left(\frac{r_k}{\max\{r_k,\sqrt{\theta}\} + \theta + B}\right) \ge \rho_{1/2} \frac{ r_k}{r_k + \theta + B} \gea \frac{\delta_k}{5.3(\delta_k + \theta + B)}\geb \frac{\delta_k}{5.3(\delta_0 + \theta + B)},  \label{eq:lb_omega1}
\end{equation}
where in (a) we use $r_k\ge \delta_k$ and $\rho_{1/2}\ge 1/5.3$ and in (b) we use the monotonicity of $\{\delta_k\}_{k\ge 0}$. 
%where in the second inequality we use %the non-decreasing property of 
%that the function $a\mapsto a/(a+\theta+B)$ is non-decreasing  on $[0,+\infty)$. 
From~\eqref{eq:decsent2}, we see that %This implies tah
\begin{equation}
\delta_{k+1}\le \bigg(1-\frac{1}{5.3(\delta_0 + \theta + B)}\bigg)\delta_k. \label{eq:lin_conv_ge}
\end{equation}
If $r_k\le \theta + B$, then $\max\{r_k,\sqrt{\theta}\}\le \theta + B$ and 
\begin{align}
\frac{r_k}{\max\{r_k,\sqrt{\theta}\} + \theta + B}\le \frac{\max\{r_k,\sqrt{\theta}\}}{\max\{r_k,\sqrt{\theta}\} + \theta + B}\le \frac{1}{2}. 
\end{align}
Therefore, by~\eqref{eq:quad_omega} in Lemma~\ref{lem:omega}, we have  
\begin{equation}
\omega\bigg(\frac{r_k}{\max\{r_k,\sqrt{\theta}\} + \theta + B}\bigg) \ge \mu_{1/2}\frac{r_k^2}{(\max\{r_k,\sqrt{\theta}\} + \theta + B)^2} \ge \frac{r_k^2}{10.6 (\theta + B)^2}, \label{eq:omega*_quad_lb}
\end{equation}
where in the last step we use $\mu_{1/2}\ge 1/2.65$ and $\max\{r_k,\sqrt{\theta}\}\le \theta + B$.

Next, if $\abst{\calS_k} = 1$, let us define $a^k := \argmax_{x\in\calS_k} \lranglet{\nabla F(x^k)}{x} = x^k$, which is consistent with the case where $\abst{\calS_k} > 1$ (cf.~Step~\ref{item:away_atom_a}). % in Algorithm~\ref{algo:AFW_SC}). 
  Note that 
\begin{equation}
r_k \ge ({1}/{2})\lranglet{\nabla F(x^k)}{a^k-v^k}. \label{eq:r_k_ave}
\end{equation}
Indeed, if $\abst{\calS_k} = 1$, then $r_k = \lranglet{\nabla F(x^k)}{a^k-v^k} = \lranglet{\nabla F(x^k)}{x^k-v^k}\ge 0$, and~\eqref{eq:r_k_ave} clearly holds; otherwise, we have $$r_k = \max\{G_k,\tilG_k\}\ge (1/2)(G_k+\tilG_k) =(1/2)\lranglet{-\nabla F(x^k)}{d_\rmF^k+d_\rmA^k} = ({1}/{2})\lranglet{\nabla F(x^k)}{a^k-v^k}.$$ 
Now, by the definition of $a^k$ and $v^k$, and Lemma~\ref{lem:lb_innerp}, we have that for any $x^*\in\calX^*$, 
\begin{equation}
 \lranglet{\nabla F(x^k)}{a^k-v^k} = {\max}_{a\in\calS_k,v\in\calV}\;\ipt{\nabla F(x^k)}{a- v}  \ge%\Phi(\calX,\calX^*) \frac{\ipt{\nabla F(x^k)}{ x^k-  x^*}}{D(\rvE x^k,\rvE x^*)}.  
\Phi(\calX,\calX^*) \frac{\ipt{\nabla F(x^k)}{ x^k-  x^*}}{\normtt{\rvE x^k - z^*}_{z^*}} \label{eq:align_lb}
\end{equation}
Combining~\eqref{eq:omega*_quad_lb},~\eqref{eq:r_k_ave} and~\eqref{eq:align_lb}, we have that for any $x^*\in\calX^*$,
\begin{align}
\omega\bigg(\frac{r_k}{\max\{r_k,\sqrt{\theta}\} + \theta + B}\bigg)\ge \frac{\lranglet{\nabla F(x^k)}{a^k-v^k}^2}{42.4 (\theta + B)^2}\ge \frac{\Phi(\calX,\calX^*)^2\ipt{\nabla F(x^k)}{ x^k-  x^*}^2}{42.4 (\theta + B)^2\normtt{\rvE x^k - z^*}_{z^*}^2}. \label{eq:omega*_lb1}
\end{align}
From the convexity of $F$, we know that $\ipt{\nabla F(x^k)}{ x^k-  x^*}\ge F(x^k)-F(x^*) = \delta_k$. Also, from Lemma~\ref{lem:eb_F}, we see that for any $x^*\in\calX^*$,
\begin{equation}
\textstyle \normtt{\rvE x^k - z^*}_{z^*}^2\le \mu_{R(y^*)}^{-1}(\barF(\rvE x^k)-\barF(z^*)) = \mu_{R(y^*)}^{-1}(F( x^k)- F(x^*)) = \mu_{R(y^*)}^{-1}\delta_k. 
\end{equation}
Plugging these into~\eqref{eq:omega*_lb1}, we have
\begin{equation}
\omega\bigg(\frac{r_k}{\max\{r_k,\sqrt{\theta}\} + \theta + B}\bigg)\ge \frac{\mu_{R(y^*)} \Phi(\calX,\calX^*)^2}{42.4 (\theta + B)^2}\delta_k. % \quad \Longrightarrow\quad  
\end{equation}
Thus, from~\eqref{eq:decsent2}, we have
\begin{equation}
\delta_{k+1}\le \bigg(1-\frac{\mu_{R(y^*)} \Phi(\calX,\calX^*)^2}{42.4 (\theta + B)^2}\bigg)\delta_k.  \label{eq:lin_conv_le}
\end{equation}
%This, together with $\ipt{\nabla F(x^k)}{ x^k-  x^*}\ge F(x^k)-F(x^*) = \delta_k$, 
%Thus, by Lemma~\ref{lem:error_bound}, %and note that $f = F$, 
Combining~\eqref{eq:lin_conv_ge} and~\eqref{eq:lin_conv_le}, in the non-full-step case, we have %$\alpha_k = b_k$, 
\begin{equation}
\delta_{k+1} \le (1-\rho) \delta_k, %\quad \mbox{where} \;\; \rho:= \min\bigg\{\frac{1}{5.3(\delta_0 + \theta + B)}, \frac{\mu_{R(y^*)} \Phi(\calX,\calX^*)^2}{42.4 (\theta + B)^2}\bigg\}< \frac{1}{5.3}. 
\end{equation}
where $\rho$ is defined in~\eqref{eq:global_lin_rate} and satisfies $\rho\le 1/5.3$ (since $\theta\ge 1$).  \\%[-1ex]

\noindent 
\underline{Combining two cases.} From the above, we know that in the full-step case, when $\alpha_k = 1$, we have $\delta_{k+1}\le (1/2)\delta_k$, and in the non-full-step case, we have $\delta_{k+1} \le (1-\rho) \delta_k$, with $\rho<1/2$.   In addition, the number of these iterations is at least $\ceil{(k-\abst{\calS_0}+q)/2}$, and in the worst case, all of them correspond to the non-full steps. This, together with the monotonicity of $\{\delta_k\}_{k\ge 0}$ in part~\ref{item:monotone}, completes the proof of part~\ref{item:lin_rate}. 
%In summary, no matter which value $\alpha_k$ takes, with at least $(k-m)/2$ iterations, we have  
%\begin{equation}
%\delta_{k+1} \le (1-\rho) \delta_k, \quad \mbox{where} \;\; \rho:= \min\left\{\frac{1}{2},\rho'\right\}. 
%\end{equation}
%end{proof}

\subsubsection{Proof of Theorem~\ref{thm:global_lin_conv_r_k} } \label{sec:proof_r_k}

%To show Theorem~\ref{thm:global_lin_conv_r_k}, we first establish a lemma that upper bounds $r_k$ in terms of $\delta_k$ (and other problem parameters). %\dom F:=\rvA^{-1}(\inter\calK)$. In addition, since $$

%The finiteness of $\bar D$ follows from the compactness of $\calS_F(x^0)$ and $\calY$

Indeed, Theorem~\ref{thm:global_lin_conv_r_k} can be viewed as a corollary of Theorem~\ref{thm:global_lin_conv}, together with a lemma that upper bounds $r_k$ in terms of $\delta_k$ (and other problem parameters), which we will establish as follows.

%The proof of Theorem~\ref{thm:global_lin_conv_G_k} 

\begin{lemma}
For all $k\ge 0$,  we have
\begin{equation}
r_k\le 4\max\big\{\delta_k,\sqrt{\delta_k}\big\}\max\{\bar D, 1\}. 
\end{equation}
\end{lemma}

\begin{proof}
Fix any $k\ge 0$ such that $G_k>0$, and hence $\delta_k> 0$. Using the same reasoning as in~\eqref{eq:curv_alpha_k}, we know that %for all $\alpha\in[0,1]$ such that $\alpha D_k<1$, 
\begin{equation}
F^*\le F(x^k+\alpha d^k)\le F(x^k) - \alpha r_k + \omega^*(\alpha D_k), \quad \forall\, \alpha\in\calI:=\{\alpha\in(0,1]:\alpha D_k<1\}. 
\end{equation}
As a result, we have 
\begin{align}
r_k\le {\inf}_{\alpha\in\calI}\; \big\{\zeta(\alpha):= {(\delta_k + \omega^*(\alpha D_k))/\alpha}\big\}. \label{eq:r_k_le_inf}
\end{align}
%To minimize $\zeta$ over $\calI$, 
Note that $\zeta'(\alpha)= (h(\alpha D_k) - \delta_k)/\alpha^2$ for $\alpha\in\calI$, where 
\begin{align}
%\zeta'(\alpha)= (h(\alpha D_k) - \delta_k)/\alpha^2,\quad \forall\,\alpha\in\calI, \;\;\; \where\;\;\; 
h(t):= 2\bar\omega^*(t) - \omega^*(t), \quad \forall\,t\in[0,1). 
\end{align}
%$\zeta'(\alpha)= (h(\alpha D_k) - \delta_k)/\alpha^2$ for $\alpha\in\calI$, where $h(t):= 2\bar\omega^*(t) - \omega^*(t)$ for $t\in[0,1)$ and $\bar\omega^*$ is defined in~\eqref{eq:omega*}. 
Note that both $h$ and $\bar\omega^*$ are strictly increasing on $[0,1)$  and nonnegative. From Lemma~\ref{lem:ub_omega*}, we have $\bar\omega^*\le h\le 2\bar\omega^*$ on $[0,1)$, and hence $(\bar\omega^*)^{-1}(\cdot/2)\le h^{-1}\le (\bar\omega^*)^{-1}$ on $\bbR_+$. 
%Additionally, 
Also, note that 
\begin{equation}
(\bar\omega^*)^{-1}(s) = \sqrt{s^2+2s} - s, \quad \forall\, s\ge 0. 
\end{equation}
%$$ for all $s\ge 0$. 
Now, let us discuss two cases. %, namely, $\delta_k> h(D_k)$ and $\delta_k\le h(D_k)$. 
%In the first case,  
If $\delta_k> h(D_k)$, then for all $\alpha\in\calI$, we have $\alpha\le 1$ and hence $\zeta'(\alpha)<0$, and from~\eqref{eq:r_k_le_inf}, we have 
\begin{equation}
r_k\le \zeta(1) = \delta_k + \omega^*(D_k)\le \delta_k + \bar\omega^*(h^{-1}(\delta_k))\le \delta_k + \bar\omega^*\big((\bar\omega^*)^{-1}(\delta_k)\big) = 2\delta_k.  \label{eq:r_k_ub0}
\end{equation}
%In the second case, 
If $\delta_k\le h(D_k)$, then $D_k>0$ and the RHS of~\eqref{eq:r_k_le_inf} has
%RHS of~\eqref{eq:r_k_le_inf} is $\zeta(\alpha^*)$ 
the optimal solution $\alpha^*:= h^{-1}(\delta_k)/D_k$. Therefore, we have
\begin{align}
%\textstyle 
r_k\le \zeta(\alpha^*) = \frac{\delta_k + \omega^*(h^{-1}(\delta_k))}{h^{-1}(\delta_k)}D_k&\lea \frac{\delta_k + \bar\omega^*((\bar\omega^*)^{-1}(\delta_k))}{(\bar\omega^*)^{-1}(\delta_k/2)}D_k= \frac{2\delta_k}{(\delta_k^2/4+\delta_k)^{1/2} - \delta_k/2}D_k, \label{eq:r_k_ub1}
%\le 2((\delta_k^2/4+\delta_k)^{1/2} + \delta_k/2), 
\end{align}
%where (a) follows from %the definition of 
%$D_k\le \bar D$ for all $k\ge 0$ (cf.~\eqref{eq:barD}). %Continuing
%From~\eqref{eq:r_k_ub1}, 
We then continue to bound $r_k$ as follows: %have
%In addition, we have 
%After algebraic manipulations, we have % it is easy to see that 
\begin{equation}
r_k \le 2((\delta_k^2/4+\delta_k)^{1/2} + \delta_k/2)D_k\le 2(\delta_k + \sqrt{\delta_k})D_k\le  4\max\big\{\delta_k,\sqrt{\delta_k}\big\}\max\{\bar D, 1\},  \label{eq:r_k_ub2}
%2\sqrt{\delta_k}(1+\sqrt{\delta_0})\bar D,
\end{equation}
where the last step follows from the monotonicity of $\{\delta_k\}_{k\ge 0}$ (cf.~Theorem~\ref{thm:global_lin_conv}\ref{item:monotone}) and~\eqref{eq:D_k_le_D}. Combining~\eqref{eq:r_k_ub0} and~\eqref{eq:r_k_ub2}, we complete the proof. 
%$\zeta'(\alpha)=0  \Longleftrightarrow  %$ if and only if $
%\alpha = h^{-1}(\delta_k)/D_k$  
%the RHS ofbecomes 
% If $\delta_k> h(D_k)$, then 
%\begin{align}
%\bar\omega^*\le h\le 2\bar\omega^* \quad \mbox{and hence}\quad (\bar\omega^*)^{-1}(\cdot/2)\le h\le (\bar\omega^*)^{-1} \quad \mbox{on $[0,1)$}. 
%\end{align}
%$$ on $[0,1)$, 
%\begin{enumerate}[label=\roman*)]
%\item $h$ is 
%\end{enumerate}
\end{proof}

\section{Local Linear Convergence}\label{sec:local_conv}

\begin{lemma} \label{lem:exist_face}
There exists a %(non-empty) 
face of $\calX$, denoted by $\calF$, such that for any  $x^*\in \calX^*$, % and any $x\in\calX$, 
if $x\in\calX$, then
%all $x\in\calX$ and all $x^*\in \calX^*$, %
%we have 
\begin{equation}
\ipt{\nabla F(x^*)}{x - x^*} = 0\quad \Longleftrightarrow \quad x\in\calF. %, \quad \forall\,x\in\calX. 
\end{equation}
\end{lemma}

\begin{proof}%\let\qed\relax
%Recall from Section~\ref{sec:char_calX*} that $\calZ = \rvE(\calX)$, $$$F(x) = \barF(\rvE x)$ and $z^*$. 
%We first show there exist
Let us first consider the problem in~\eqref{eq:poi_yt}, namely $\min_{z\in\calZ}\, \barF(z)$, which has a unique solution $z^*\in\calZ$ (cf.~Section~\ref{sec:char_calX*}). From the first-order optimality condition, we know that $-\nabla \barF(z^*)\in\calN_{\calZ}(z^*)$, where $\calN_{\calZ}(z^*):= \{d\in\bbY\times\bbR:\ipt{d}{z-z^*}\le 0\}$ denotes the normal cone of $\calZ$ at $z^*$. Let 
$$\calZ=\{z\in\bbY\times\bbR: \ipt{c_i}{z}\le h_i, \;\forall\,i\in[m]\}$$ 
be the H-representation of $\calZ$, and let $\calI^*:=\{i\in[m]: {c_i}{z^*}= h_i\}\subseteq [m]$ be the indices of the active constraints at $z^*$. From standard results (see e.g., ),  we know that $\calN_{\calZ}(z^*):= \cone\{c_i\}_{i\in\calI^*}$, and hence $-\nabla \barF(z^*) = \sum_{i\in\calI^*}\, \gamma_i c_i$ for some $\{\gamma_i\}_{i\in\calI^*}\subseteq \bbR_+$. Define $\calI:=\{i\in\calI^*:\gamma_i>0\}\subseteq \calI^*$, and 
\begin{align}
\ipt{\nabla F(z^*)}{z - z^*} = \textstyle\sum_{i\in\calI}\, \gamma_i \ipt{c_i}{z^*-z}= \textstyle\sum_{i\in\calI}\, \gamma_i (h_i-\ipt{c_i}{z}).
\end{align}
For any $z\in\calZ$, if $\ipt{\nabla F(z^*)}{z - z^*}=0$, then $\ipt{c_i}{z} = h_i$ for all $i\in\calI$, and the converse is clearly true. Define 
\begin{equation}
\calF_\calI:=\{z\in\bbY\times\bbR: \ipt{c_i}{z}= h_i, \;\forall\,i\in\calI, \;\; \ipt{c_i}{z}\le h_i, \;\forall\,i\in[m]\setminus\calI\}, \label{eq:calF}
\end{equation}
 which is a face of $\calZ$. From above, we know that for all $z\in\calZ$, 
\begin{equation}
\ipt{\nabla \barF(z^*)}{z - z^*} = 0\quad \Longleftrightarrow \quad z\in\calF_\calI. \label{eq:equiv_Z}
\end{equation}
%From Section~\ref{sec:char_calX*}, 
Recall that $\calZ:= \rvE(\calX)$ and $\{z^*\} = \rvE(\calX^*)$. 
In Lemma~\ref{lem:exist_face}, we can  take $\calF:= \{x\in\calX:\rvE x \in\calF_\calI\}$, namely the pre-image of $\calF_\calI$ under $\rvE$ w.r.t.\ $\calX$, and note that $\calF$ is indeed a face of $\calX$. 
From~\eqref{eq:equiv_Z}, we know that for all $x\in\calX$ and   all $x^*\in\calX^*$, \phantom{\qedhere}
\begin{align}
x\in\calF \;\;\Longleftrightarrow\;\; \rvE x \in \calF_\calI \;\;\Longleftrightarrow\;\; \ipt{\nabla \barF(\rvE x^*)}{\rvE x - \rvE x^*} = 0 \;\;\Longleftrightarrow\;\;\ipt{\nabla F(x^*)}{x - x^*} = 0. \tag*{$\qed$}
\end{align}

\end{proof}

\begin{remark} \label{rmk:minimal_face}
Clearly, we have $\calX^*\subseteq \calF$. In addition, note that %the (affine) dimension of the face 
$\calF$ can be potentially much smaller than $\calX$. In particular, if strict complementarity holds for the problem in~\eqref{eq:poi_yt}, i.e., $-\nabla \barF(z^*)\in\ri \calN_{\calZ}(z^*)$, then $\calF$ (as constructed in the proof of Lemma~\ref{lem:exist_face}) is the minimal face of $\calX$ that contains $\calX^*$. For a proof of this statement, see Appendix~\ref{app:details_rmk}. 
\end{remark}

Based on Lemma~\ref{lem:exist_face}, define $\calV_\calF:= \calV\cap\calF$, $\bar\calV_\calF:= \calV\setminus\calF$, $\Delta_\calF(x^*):= {\min}_{v\in\bar\calV_\calF}\, \ipt{\nabla F(x^*)}{v - x^*}$ for $x^*\in\calX^*$, %and 
\begin{equation}
\hspace{-.6ex}
\barcalX^*:= {\argmax}_{x^*\in\calX^*} \Delta_\calF(x^*)\quad \mbox{and}\quad\barDelta_\calF := %{\max}_{x^*\in\calX^*}\;
\Delta_\calF(x^*) \;\;\mbox{for}\;\;x^*\in\barcalX^*,
%\Delta_\calF:= {\min}_{v\in\bar\calV_\calF}\; \ipt{\nabla F(x^*)}{v - x^*}>0. 
 \label{eq:def_Delta_F}
\end{equation}
where we know that $\barcalX^*\ne \emptyset$, since the function $\Delta_\calF$ is continuous on the compact set $\calX^*$. %the existence of $\barx^*\in\calX^*$ follows from the continuity of $\Delta_\calF$ on the compact set $\calX^*$. %, which is nonempty and compact. %the compactness of $\calX^*$ and . 

%\noindent 
%The next lemma 
Next, we show that if $x^k$ is ``close'' to $\calX^*$, in the sense that %the objective gap 
$\delta_k$ is small, then for all $x^*\in\calX^*$, the functions $x\mapsto \ipt{\nabla F(x^k)}{x-x^k}$ and $x\mapsto \ipt{\nabla F(x^*)}{x-x^*}$ are ``close'' {\em uniformly} on $\calX$.

\begin{lemma}\label{lem:delta_k_ub}
Fix any $\Delta>0$. In Algorithm~\ref{algo:AFW_SC}, if 
\begin{equation}
\delta_k<  V(\Delta):= \min\bigg\{\frac{\sqrt{\mu_{R(y^*)}}\Delta}{2R(y^*)+\Delta},\;\frac{\Delta R(y^*)}{4((2R(y^*)+\Delta)R(y^*)+\Delta)}\bigg\}^2, % < \frac{1}{16},
\label{eq:delta_k_le}
\end{equation}
then for all $x^*\in\calX^*$, we have 
\begin{equation}
%{\sup}_{x^*\in\calX^*}\;
{\sup}_{x\in\calX}\;\abst{\ipt{\nabla F(x^k)}{x-x^k} - \ipt{\nabla F(x^*)}{x-x^*}}<\Delta. \label{eq:unif_ipt_bound}
\end{equation}
%$$. 
%and $\mu_{R(y^*)}$ are defined in~\eqref{eq:def_R_y*} and~\eqref{eq:quad_omega}, respectively. 
%$\delta_$
\end{lemma}
\begin{proof}
See Appendix~\ref{app:proof_delta_k_ub}. 
\end{proof}

\begin{remark}
Several remarks are in order. 
First, note that the upper bound $V(\Delta)$ in~\eqref{eq:delta_k_le} only depends on $R(y^*)$, namely the radius of $\calY$ centered at $y^*$ under $\normt{\cdot}_{y^*}$. Second, there might exist several different expressions for  $V(\Delta)$, and the one in~\eqref{eq:delta_k_le} simply provides an example. Third, note that if the linear function $\ipt{c}{\cdot}$ is constant on $\calX$ (in which case we have $\ipt{\nabla F(x)}{x}= -\theta$ for all $x\in\calX$), then the expression of $V(\Delta)$ can be made much simpler. This happens, for example, in D-optimal design (cf.~\eqref{eq:Dopt} and~\eqref{eq:Dopt2}) and inference of MHP (cf.~\eqref{eq:L1_MLE_rewrite}).  
\end{remark}

Now, let us present the main theorem of this section, which states that when $\delta_k$ falls below a threshold, the set of active atoms $\calS_k$ forms a {\em nested downward} sequence, until $x^k$ lands on $\calF$. In addition, once $x^k\in\calF$, it will remain in $\calF$.  

\begin{theorem}\label{thm:local_lin_conv}
Define $\delta_{\min}(\bar\calV_\calF):= \min_{v\in\bar\calV_\calF} \delta(v)>0$, where $\delta(\cdot)$ is defined in~\eqref{eq:delta}. 
In Algorithm~\ref{algo:AFW_SC}, let $\bar k\ge 0$ satisfy that $\delta_{\bar k}< \min\{ V(\barDelta_\calF/2), \delta_{\min}(\bar\calV_\calF)\}$. % for some , %where $\delta_{\min}(\bar\calV_\calF):= \min_{v\in\bar\calV_\calF} \delta(v)$, 
For all $k\ge \bar k$,  if $x^k\not\in\calF$, then % we have 
\begin{enumerate}[label=(\roman*)]
\item \label{item:nested} $\calS_{k+1}\subseteq\calS_{k}$, when either exact line-search or adaptive stepsize is used in Step~\ref{item:step_size_SC_a},
\item \label{item:nested_strict} $\calS_{k+1}=\calS_{k}\setminus\{a^k\}$ for some $a^k\in \calS_{k}\cap\bar\calV_\calF$, when exact line-search is used in Step~\ref{item:step_size_SC_a}; 
\end{enumerate}
otherwise, if $x^k\in\calF$, then $x^{l}\in \calF$ for all $l\ge k$. 
\end{theorem}

\begin{proof}
Since $x^k\not\in\calF$, we have $\calS_k\cap \bar\calV_\calF\ne \emptyset$. 
By the monotonicity of $\{\delta_k\}_{k\ge 0}$, we have $\delta_k< V(\barDelta_\calF/2)$ for all $k\ge \bar k$. In Lemma~\ref{lem:delta_k_ub}, take any $x^*\in\barcalX^*$,  and by the definition of $\barDelta_\calF$ in~\eqref{eq:def_Delta_F}, we have 
\begin{align}
{\min}_{v\in\bar\calV_\calF}\;\lranglet{\nabla F (x^k)}{v-x^k}> \barDelta_\calF/2  &> {\max}_{v\in\calV_\calF}\;\lranglet{\nabla F (x^k)}{x^k-v}\label{eq:for_any_x0}\\
&\ge {\min}_{v\in\calV_\calF}\;\lranglet{\nabla F (x^k)}{x^k-v}> -\barDelta_\calF/2. \label{eq:for_any_x}
\end{align}
As a result, we have ${\max}_{v\in\calV_\calF}\;\lranglet{\nabla F (x^k)}{v-x^k}< \barDelta_\calF/2$, and hence i) if $d^k=d_\rmA^k$, % (i.e., away step is selected), 
then $a^k\in \bar\calV_\calF$ and ii) $\tilG_k= {\max}_{v\in\calS_k\cap\bar\calV_\calF}\;\lranglet{\nabla F (x^k)}{v-x^k} > \barDelta_\calF/2$. In addition, since ${\max}_{v\in\bar\calV_\calF}\;\lranglet{\nabla F (x^k)}{x^k-v}< -\barDelta_\calF/2$, we have $G_k = \max_{v\in\calV_\calF} \lranglet{\nabla F (x^k)}{x^k-v}< \barDelta_\calF/2<\tilG_k$. Since $\delta_k< \delta_{\min}(\bar\calV_\calF)$, we know that $x^k\not\in\bar\calV_\calF$, but since $\calS_k\cap \bar\calV_\calF\ne \emptyset$, we have $\abst{\calS_k}>1$. Hence, in Step~\ref{item:choose_direction}, we choose $d^k=d^k_\rmA$, which implies that $\calS_{k+1}\subseteq\calS_{k}$. Also, from the discussion above, we have $a^k \in\bar\calV_\calF$. This proves~\ref{item:nested}. %. Thus, we have 
%Note that this implies that $\calS_{k+1}\subseteq\calS_{k}$. %, which holds for both  exact line-search and adaptive stepsize. 

We can strengthen the result above if exact line-search is used, in which case we have $\alpha_k = \bar\alpha_k$, and $a^k$ is dropped from $\calS_k$. To see this, suppose that $\alpha_k < \bar\alpha_k$, then we have 
\begin{align}
(1+\alpha_k)^{-1}\ipt{\nabla F(x^{k+1})}{a^k - x^{k+1}} = \ipt{\nabla F(x^{k+1})}{a^k - x^k}=0. \label{eq:ipt_x^k+1}
\end{align}
Since $\delta_{k+1}\le \delta_k<V(\barDelta_\calF/2)$, we have ${\min}_{v\in\bar\calV_\calF}\;\lranglet{\nabla F (x^{k+1})}{v-x^{k+1}}> \barDelta_\calF/2$, and since $a_k \in\bar\calV_\calF$, this contradicts~\eqref{eq:ipt_x^k+1}. This proves~\ref{item:nested_strict}. 

Finally, if $x^k\in\calF$ and $d^k=d^k_\rmF$, since ${\max}_{v\in\bar\calV_\calF}\;\lranglet{\nabla F (x^k)}{x^k-v}< -\barDelta_\calF/2$, we must have $v^k\in\calV_\calF$. This shows $x^{k+1}\in\calF$, and we complete the proof. 
%$$$$
%\begin{align*}
%\tilG_k= {\max}_{v\in\bar\calV_\calF}\;\lranglet{\nabla F (x^k)}{v-x^k} \ge 2\Delta_\calF/3.  
%\end{align*}
%and ii)   
%\begin{align*}
%{\max}_{v\in\calV}\;\lranglet{\nabla F (x^k)}{v-x^k} = \max\{{\max}_{v\in\calV_\calF}\;\lranglet{\nabla F (x^k)}{v-x^k},{\max}_{v\in\bar \calV_\calF}\;\lranglet{\nabla F (x^k)}{v-x^k}\}
%\end{align*}
%Note that by taking $\Delta = \Delta_\calF/3$ in~\eqref{eq:unif_ipt_bound}, we have
\end{proof}

\begin{remark} \label{rmk:better_local_rate}
Theorem~\ref{thm:local_lin_conv} implies that if exact line-search is used, then $\{x^k\}_{k\ge K}\subseteq\calF$ for $K= \bar k + \abst{\calS_{\bar k}\cap\bar\calV_\calF}$. As a result, we can apply Theorem~\ref{thm:global_lin_conv} to analyze  the linear convergence rate of $\{\delta_k\}_{k\ge K}$ by treating $\calF$ as $\calX$. %letting $\calX=\calF$. 
Note that this local linear rate can be potentially much better than the global one, 
%the tail sequence $\{x^k\}_{k\ge K}$ and the linear convergence rate of the tail sequence $\{x^k\}_{k\ge K}$ on $\calF$ 
%for some $K\ge 0$, we have .  
since $\calF$ can be potentially much smaller than $\calX$ (cf.\ Remark~\ref{rmk:minimal_face}), and %on $\calF$, %this may potentially improve 
several quantities involved in the definition of the linear rate   $\rho$ (e.g., $\Phi(\calX,\calX^*)$,  $R(y^*)$ and $B$) may potentially become better when defined on $\calF$ (cf.~\eqref{eq:global_lin_rate}). 
\end{remark}

\section{Numerical Experiments} \label{sec:experiments}

In this section, we investigate the numerical performance of our away-step FW method (i.e., Algorithm~\ref{algo:AFW_SC}) on two  applications as described in Section~\ref{sec:intro}, namely, the D-optimal design problem~\eqref{eq:Dopt} and the (reformulated) MHP inference problem~\eqref{eq:L1_MLE_rewrite}. All of the experiments were conducted in Python 3.8 on a machine with an Intel i7-9700 3-GHz CPU and 12-GB RAM. 

\subsection{Experimental Setup}\label{sec:setup}

{\bf Algorithms under comparison.} As mentioned in Section~\ref{sec:intro},  both~\eqref{eq:Dopt} and~\eqref{eq:L1_MLE_rewrite} are  rather ``non-standard'' from the FOM point of view %. Specifically,  
--- in either problem, the objective function $F$ does not have Lipschitz gradient on the constraint set $\calX$. Due to this reason, these two problems have been rarely studied in the FOM literature. That said, as reviewed in Section~\ref{sec:lit}, some principled FOMs (with convergence guarantees) have been recently proposed and are able %that be applied to these 
to solve these two problems. These algorithms are used for comparison with Algorithm~\ref{algo:AFW_SC} and are listed below:
%For both applications, we shall 

\begin{itemize}[itemsep=0ex]
	\item {\tt FW-E} -- the FW method for minimizing $\theta$-LHSCB with the exact line-search~\cite{Dvu_20,Zhao_23},
	\item {\tt FW-A} -- the FW method for minimizing $\theta$-LHSCB with the adaptive step-size~\cite{Dvu_20,Zhao_23},
	\item {\tt RSGM-F} -- relatively smooth gradient method with fixed step-size~\cite{Bauschke_17,Lu_18}, 
	\item {\tt RSGM-B} -- relatively smooth gradient method with backtracking line-search~\cite{Ston_21},
	\item {\tt MG} -- the multiplicative gradient method~\cite{Ari_72,Silvey_78,Zhao_22mg}.\footnote{Note that the multiplicative gradient method has long been proposed in the statistics community for solving~\eqref{eq:Dopt} and~\eqref{eq:L1_MLE_rewrite} separately~\cite{Ari_72,Silvey_78}, and Zhao~\cite{Zhao_22mg} recently generalized this method to a class of problems.} 
%	a simple algorithm developed by Cover in 1984 specifically for a problem that is equivalent to the normalized PET problem ~\cite{Cover_84}.  %\red{[Why do you call it EM?]}
\end{itemize}
%We excluded mirror descent from our computational comparisons because the sparsity of $P$ violates the basic assumption needed to apply the mirror descent (MD) method to \eqref{eq:PET_final} (see e.g.,~\cite{BenTal_01}).  We now review the relevant details of the three algorithms {\tt RSGM-F}, {\tt RSGM-LS} and {\tt EM}. 
In addition, for ease of reference, we list two variants of Algorithm~\ref{algo:AFW_SC} below:
\begin{itemize}
	\item {\tt AFW-E} -- Algorithm~\ref{algo:AFW_SC} with the exact line-search (cf.\ Step~\ref{item:step_size_SC_a}\ref{item:exact_LS}),
	\item {\tt AFW-A} -- Algorithm~\ref{algo:AFW_SC} with the adaptive step-size (cf.\ Step~\ref{item:step_size_SC_a}\ref{item:adapt_step}).
\end{itemize}	
Note that all of these methods above are {\em feasible}, namely, when starting from a feasible point in $\calX\cap\dom F$,  %each method will generate a
the sequence of iterates  generated by each method will remain in $\calX\cap\dom F$. 
In addition, for both~\eqref{eq:Dopt} and~\eqref{eq:L1_MLE_rewrite}, exact line-search (cf.\ Step~\ref{item:step_size_SC_a}\ref{item:exact_LS}) can be efficiently implemented. For details, see~\cite{Zhao_23} and~\cite{Khachiyan_96}. % with low computational cost. 
\\

%since in Step~\ref{item:step_size_SC_a} of Algorithm~\ref{algo:AFW_SC}, 
%similar to the FW method in~\cite{Zhao_23}, the step-size in Algorithm~\ref{algo:AFW_SC}, % $\alpha_k$
% can be chosen via two ways, and we refer to Algorithm~\ref{algo:AFW_SC} with adaptive step-size and exact line-search 

\noindent {\bf Metrics of comparison.} We will primarily use the {\em objective gap} $\{\delta_k\}_{k\ge 0}$ %(i.e., $F(x) - F^*$ for $x\in\calX$) 
to compare the convergence speed of the methods mentioned above, where recall that $\delta_k:= F(x^k) - F^*$.  Since the definition of $\delta_k$ %\, (:= F(x^k) - F^*)$ 
involves  $F^*$, which is typically unknown, we first run {\tt AFW-E} (or in principle any FW-type algorithm) to obtain an iterate $x^k\in\calX$ with FW gap $G_k < \epsilon$, 
which then implies that  $F(x^k) <  F^*+\epsilon$ (cf.~\eqref{eq:def_G_k}). By choosing $\epsilon$ sufficiently small, $F(x^k)$ then provides a reasonably good approximation of $F^*$. For the purpose of our experiments, we choose $\epsilon:=  10^{-9}$. %Besides the objective gap, 

In addition, since %the constraint set 
$\calX$ is the unit simplex 
in both~\eqref{eq:Dopt} and~\eqref{eq:L1_MLE_rewrite}, we are also interested in comparing the {\em sparsity} (i.e., the number of non-zeros) of the iterates $\{x^k\}_{k\ge 0}$ generated by each method. Indeed, for both problems, a sparse (approximately optimal) solution is desired. For example, when~\eqref{eq:L1_MLE_rewrite} is used to analyze a social network, a sparse solution $a$ %to~\eqref{eq:L1_MLE_rewrite}   
helps reveal the underlying sparse network structure~\cite{Zhou_13}. % 
Note that in our experiments, to accommodate finite machine precision, we define an entry of $x^k$ to be non-zero if it is greater than the machine accuracy. 

Finally, for FW-type methods (including {\tt FW-E}, {\tt FW-A}, {\tt AFW-E} and {\tt AFW-A}), we will also compare the FW gap $\{G_k\}_{k\ge 0}$ (cf.~\eqref{eq:G_k}) generated by these methods.

\subsection{First Experiment: D-Optimal Design} %~\eqref{eq:Dopt} }

In~\eqref{eq:Dopt}, we set $m=2000$, $n=100$ and generated the $m$ %problem data 
vectors $\{a_i\}_{i=1}^m \subseteq \mathbb{R}^n$ independently such that  for each $i\in[m]$, $a_i\sim\calN(0,10I_n)$, namely the centered normal distribution with covariance matrix $10I_n$. We compare the numerical performance of all the methods listed in Section~\ref{sec:setup} with starting point $x^0 = (1/m)e$ (recall that $e:= (1,\ldots,1)$), namely the analytic center of $\Delta_m$. %For all methods, 
% and $a_i$'s are independent. 
The results are shown in Figure~\ref{fig:Dopt}, from which we can make the following observations. First, from Figure~\ref{fig:Dopt_obj}, we see that in terms of objective gap,  both {\tt AFW-E} and {\tt AFW-A} converge linearly, and are significantly faster than the other methods, which converge sub-linearly. In addition, the linear rates of both methods are improved somewhere between 1s and 2s.  Indeed, this can be precisely explained examining Figure~\ref{fig:Dopt_sp}, where we observe that around the same time, both methods land on (and remain in) some low-dimensional faces of $\Delta_m$, and according to Remark~\ref{rmk:better_local_rate}, this leads to an improved local linear rate. Second, from Figure~\ref{fig:Dopt_sp}, it is evident that {\tt AFW-E} and {\tt AFW-A} can return much sparser solutions compared to the other methods. In fact, for all of the other methods except {\tt MG}, the sequence of iterates $\{x^k\}_{k\ge 0}$ generated are fully dense (i.e., they remain in $\ri\Delta_m$). Third, from Figure~\ref{fig:Dopt_FWG}, in terms of the FW gap $\{G_k\}_{k\ge 0}$, {\tt AFW-E} and {\tt AFW-A} also converge linearly, and are much faster compared to {\tt FW-E} and {\tt FW-A}. In addition, for both {\tt AFW-E} and {\tt AFW-A}, the slope of the (straight line) plot in Figure~\ref{fig:Dopt_FWG} is about 1/2 of that in Figure~\ref{fig:Dopt_obj}, which precisely corroborates our theory --- see Theorem~\ref{thm:global_lin_conv_r_k} and Remark~\ref{rmk:1/2-worse}. 
%, which states 
%by comparing the slope of the 
 %Indeed, both of them reach an accuracy of $10^{-8}$ in about six seconds. 

\begin{figure}[t]
     \centering
     \begin{subfigure}[b]{0.333\columnwidth}
         \centering
         \includegraphics[width=\textwidth]{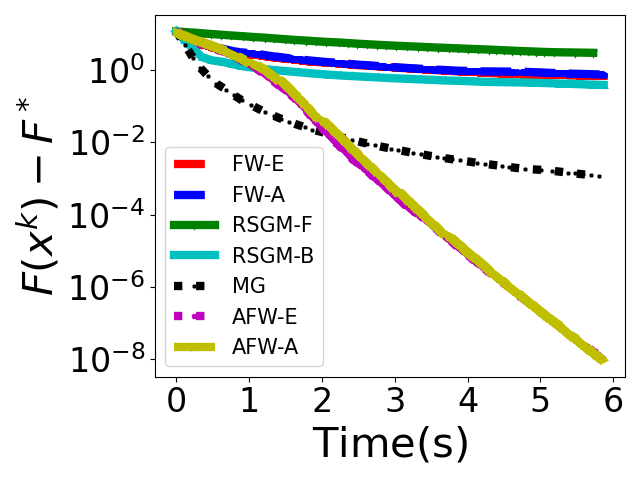}
         \caption{Log-linear plot of obj.\ gap}\label{fig:Dopt_obj}
     \end{subfigure}\hfill
     \begin{subfigure}[b]{0.333\columnwidth}
         \centering
         \includegraphics[width=\textwidth]{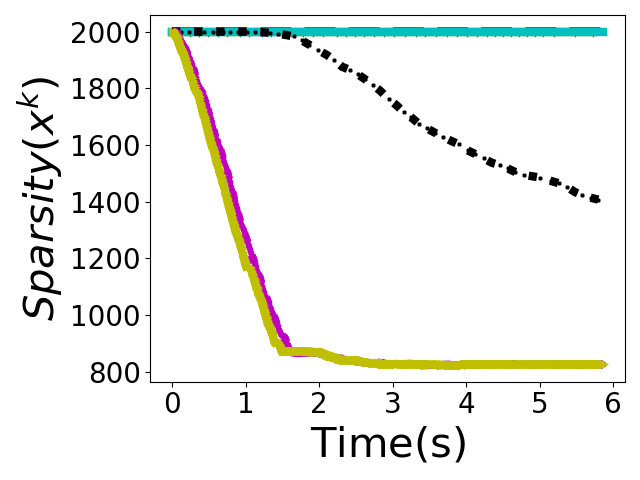}
         \caption{Plot of sparsity}
         \label{fig:Dopt_sp}
     \end{subfigure}\hfill
     \begin{subfigure}[b]{0.333\columnwidth}
         \centering
         \includegraphics[width=\textwidth]{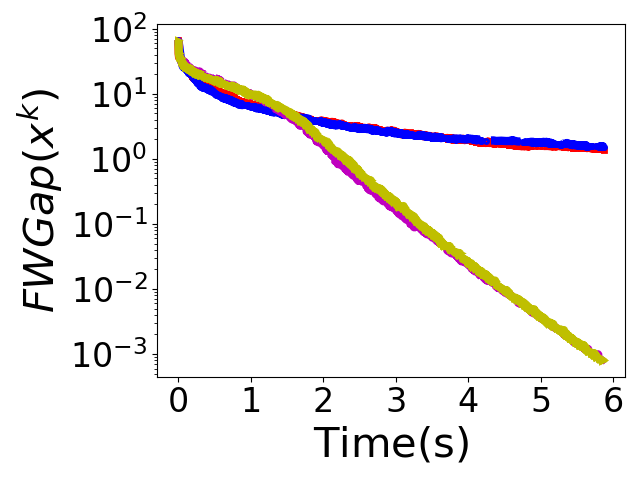}
         \caption{Log-linear plot of FW gap}
         \label{fig:Dopt_FWG}
     \end{subfigure}
        \caption{Results of the D-optimal design problem~\eqref{eq:Dopt}.}
        \label{fig:Dopt}
\end{figure}

\subsection{Second Experiment: Inference of MHP}

We first generated the arrival points of a $m$-dimensional MHP over time interval $[0,t)$ with parameters $\bar\mu$ and $A$ 
 %$\{\mu_k\}_{k\in[m]}$ and $\{a_{k,l}\}_{k,l\in[m]}$ 
(cf.\ Section~\ref{sec:two_other}) using the simulator in~\cite{Morse_20}. We set $m=1000$, $t = 5000$ and $\bar\mu = 0.1 e$. %the base intensities $\mu_k = 0.1$ for all $k\in[m]$. 
To generate a sparse infectivity matrix $A$, we first drew %$\hat A\in\bbR^{m\times m}$ with 
its entries independently from %${\sf Unif}[0.1, 0.5]$
the uniform distribution over $[0.1, 0.5]$, and then uniformly randomly set 90\% of its entries to zero. For stability of simulation, we also normalized $A$ so that its maximum eigenvalue equaled 0.9.  By denoting $x:= (\mu,a)$ in~\eqref{eq:L1_MLE_rewrite}, we compare the numerical performance of all the methods listed in Section~\ref{sec:setup} with starting point $x^0 := (m+1)^{-1}e$, and the results are shown in Figure~\ref{fig:MHP}. Indeed, we can make the same observations from Figure~\ref{fig:MHP} as those from Figure~\ref{fig:Dopt}. In particular, compared to the other methods, {\tt AFW-E} and {\tt AFW-A} not only are  more efficient in reducing the objective gap (and the FW gap), but also are more effective in producing a sparse solution. In addition, empirical behaviors of {\tt AFW-E} and {\tt AFW-A} corroborate our theoretical results (cf.\ Sections~\ref{sec:global_lin_conv} and~\ref{sec:local_conv}) fairly well. 
%  are fairly consistent with 
%We generated the entries of the infectivity matrix $A$ independently  by 

\begin{figure}[t]
     \centering
     \begin{subfigure}[b]{0.333\columnwidth}
         \centering
         \includegraphics[width=\textwidth]{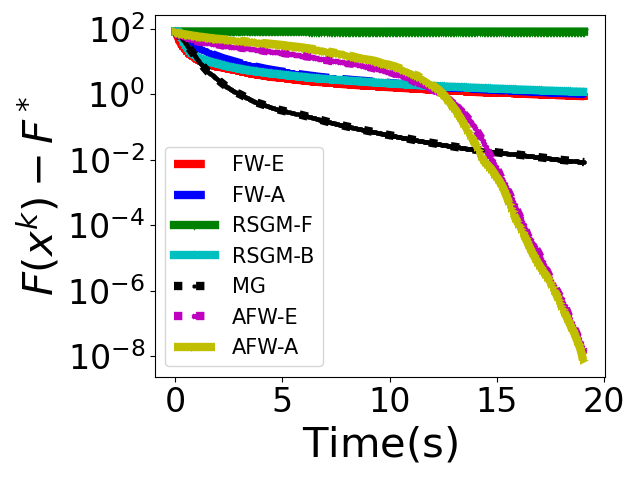}
         \caption{Log-linear plot of obj.\ gap}
         \label{fig:MHP_obj}
     \end{subfigure}\hfill
     \begin{subfigure}[b]{0.333\columnwidth}
         \centering
         \includegraphics[width=\textwidth]{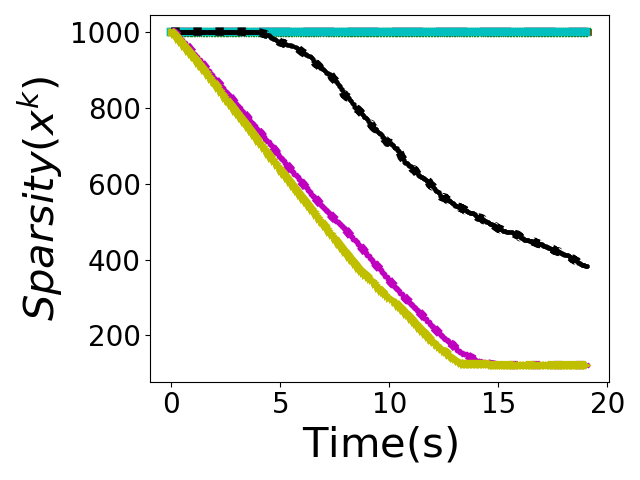}
         \caption{Plot of sparsity}
         \label{fig:MHP_sp}
     \end{subfigure}\hfill
     \begin{subfigure}[b]{0.333\columnwidth}
         \centering
         \includegraphics[width=\textwidth]{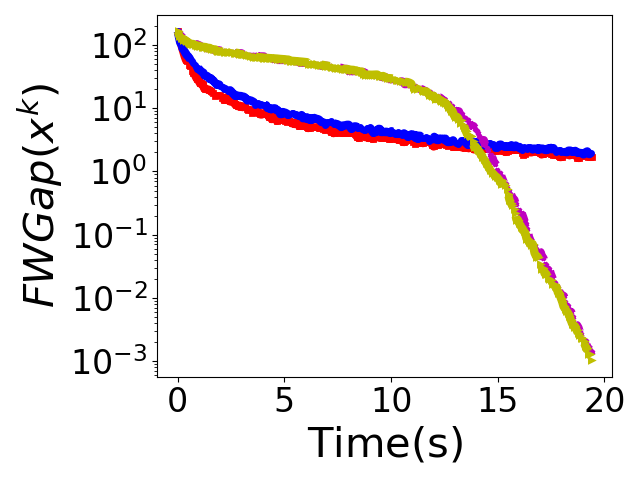}
         \caption{Log-linear plot of FW gap}
         \label{fig:MHP_FWG}
     \end{subfigure}
        \caption{Results of the MHP inference problem~\eqref{eq:L1_MLE_rewrite}.}
        \label{fig:MHP}
\end{figure}

%Following the description of a 
%
%We set $m=1000$ and $t = 5000$, and generated the arrival points of a $m$-dimensional MHP over time interval $[0,t)$ with in Section~\ref{sec:two_other},  and 

%\begin{figure}
%\includegraphics[width=.333\columnwidth]{m=2000_n=100_obj_time.png} %\hspace{-1em}
%\includegraphics[width=.333\columnwidth]{m=2000_n=100_sp_time.png} 
%\includegraphics[width=.333\columnwidth]{m=2000_n=100_FWG_time.png} 
%\end{figure}

%or $x\in\Delta_n$, we define its number of non-zeros as 
%, and we plot $\{\delta_k\}_{k\ge 0}$ until it reaches $10\epsilon$
%approximate optimal solution $\hatx\in\calX$ such that $G(\hatx)<  10^{-9}$. 

%\newpage

\section*{Acknowledgment}

The author is grateful for the help discussions with Prof.\ Robert M.\ Freund during the preparation of this manuscript. The author's research is supported by AFOSR Grant No.\ FA9550-22-1-0356.

%\newpage

\appendix

\section{Proof of Theorem~\ref{lem:omega} }\label{app:proof_omega}

First, let us show that 
\begin{equation}
\omega(t)\ge \frac{t^2}{2(1+t)}, \quad\forall\,t\ge 0. 
\end{equation}
%$t\ge 0$. 
To see this, note that $\omega(0)=0$ and hence for all $t\ge 0$, we have
\begin{equation}
\omega(t) = \int_{0}^t \; \omega'(s)\; \rmd s = \int_{0}^t \; \frac{s}{1+s}\; \rmd s\ge \int_{0}^t \; \frac{s}{1+t}\; \rmd s  = \frac{t^2}{2(1+t)}. \label{eq:lb_omega}
\end{equation}
To show~\eqref{eq:quad_omega}, %it suffices to 
we first show that the function $\zeta(t):= \omega(t)/t^2$ is non-increasing on $(0,+\infty)$. Indeed,  from~\eqref{eq:lb_omega}, we  have that for all $t> 0$, 
\begin{equation}
\zeta'(t) = \frac{1}{2 t^3}\bigg(\frac{t^2}{2(1+t)} - \omega(t)\bigg)\le 0. 
\end{equation}
Therefore, for all $t\in(0,\beta]$, we have $\omega(t)\ge (\omega(\beta)/\beta^2)  t^2 = \mu_\beta t^2$, and this also holds at $t=0$ since $\omega(0)=0$. % --- this shows~\eqref{eq:quad_omega}. 
%$\zeta'(t) = \xi(t)/t^3$, where $\xi(t):= t^2/(1+t) - 2\omega(t)\le 0$. 
Similarly, to show~\eqref{eq:lin_omega}, it suffices to show that  the function $\bar\zeta(t):= \omega(t)/t$ is non-decreasing on $(0,+\infty)$. Indeed, for all $t> 0$, from~\eqref{eq:lb_omega}, we  have that for all $t> 0$, 
\begin{equation}
\bar\zeta'(t) = \frac{1}{t^2}\bigg(\ln(1+t) - \frac{t}{1+t}\bigg)\ge \frac{1}{t^2}\bigg(1-\frac{1}{1+t} - \frac{t}{1+t}\bigg)= 0, 
\end{equation}
where we use $\ln(x)\ge 1-1/x$ for all $x>0$.

\section{Proof of Lemma~\ref{lem:local_norm_grad_diff} } \label{app:proof_grad_diff}

For any $\bary\in\calB_y(y,1)$, let us define the local norm of $\nabla^2 f(\bary)$ at $y$ as
\begin{align}
\normt{\nabla^2 f(\bary)}_{y}:= {\sup}_{\normt{u}_y=1}\; \normt{\nabla^2 f(\bary) u}_{y}^* = {\sup}_{\normt{u}_y=1}\; {\ipt{\nabla^2 f(\bary) u}{u}}= {\sup}_{u\ne 0}\; \normt{u}_{\bary}^2/\normt{u}_y^2. 
\end{align}
From Lemma~\ref{lem:approx_Hessian}, we know that $\normt{\nabla^2 f(\bary)}_{y}\le (1-\normt{\bary-y}_y)^{-2}$. Therefore, we have \phantom{\qedhere}
\begin{align}
\normt{\nabla f(y') - \nabla f(y)}_{y}^* &\le \int_{0}^1\; \normt{\nabla^2 f(y+t(y'-y))}_y\normt{y'-y}_y \rmd t \nn\\
&\le \int_{0}^1\; \frac{\normt{y'-y}_y}{(1-t\normt{y'-y}_y)^{2}} \rmd t = \frac{\normt{y'-y}_y}{1-\normt{y'-y}_y}. \tag*{$\qed$}
\end{align}

\section{Proof of Lemma~\ref{lem:lb_Phi}} % $\Phi(\calX,\calX^*)>0$} 
\label{app:positivity_Phi}

We will need the following elementary lemmas in our proof. For completeness, we provide their proofs at the end of this section. 

\begin{lemma}\label{lem:conv_face}
Let $\calV\subseteq \bbY$ be a nonempty  finite set %(i.e., $1\le \abst{\calF}<+\infty$) 
and let $\calF$ be any face of %\in{\sf Faces}
$\conv(\calV)$. 
Then we have $\conv(\calV\cap \calF)= \calF$. 
\end{lemma}

\begin{lemma}\label{lem:conv_argmin}
Let $\calV\subseteq \bbY$ be a nonempty  finite set.  %(i.e., $1\le \abst{\calF}<+\infty$) 
Then for any $p\in\bbY$, we have %$\conv(\calV\cap \calF)= \calF$. 
\begin{align}
\textstyle \conv(\argmin_{v\in\calV}\, \ipt{p}{v}) = \argmin_{x\in\conv(\calV)}\, \ipt{p}{v}. 
\end{align}
\end{lemma}

%\begin{proof}
%
%\end{proof}

Let $x\in\calX\setminus\calX^*$ and $\calS\in {\sf AR}(x)$. For convenience, define $\bar\calS:= \rvE(\calS)$, $z\ne z^*$ and $\bar\calV:= \rvE(\calV)$, and we have $z:=\rvE x\in\conv(\bar\calS)$ and $\calZ= \conv(\bar\calV)$. 
To begin with, let us note that 
\begin{alignat}{2}
\Phi_\calS(x,\calX^*) = &\;{\min}_{p:\ipt{p}{d}=1}\; &&{\max}_{a\in\bar\calS,v\in\bar\calV}\; \ipt{p}{a-  v} \\
\ge &\;{\min}_{p:\ipt{p}{d}=1}\; &&{\max}_{a\in\bar\calS,v\in\bar\calV\,\cap \,\calF_{z^*}}\; \ipt{p}{a-  v}  \label{eq:Phi_S_local}
\end{alignat}
Note that the RHS of~\eqref{eq:Phi_S_local} can be equivalently written as the following LP: 
\begin{equation}
\begin{split}
{\min}_{p,t,\tau} \;\;\;\;&  t-\tau  \\
\st\;\;\;\;  &\ipt{p}{a} \le t, \quad\;\;\; \forall\, a\in\bar\calS\\
& \ipt{p}{v}\ge \tau, \quad\;\;\; \forall\, v\in\bar\calV\cap \calF_{z^*}\\
& \lranglet{p}{d} = 1
\end{split} \tag{PLP}\label{eq:PLP}
\end{equation}
The dual of~\eqref{eq:PLP} can reads
\begin{equation}
\begin{split}
{\max}_{\lambda,\pi,\gamma} \;\; &\gamma\\
\st \;\; &\textstyle \sum_{a\in\bar\calS} \lambda_a = 1, \;\;\;\;\qquad \lambda_a\ge 0, \quad \forall\, a\in \bar\calS,\\
&\textstyle \sum_{v\in\bar\calV \,\cap \,\calF_{z^*}} \pi_v = 1, \quad  \pi_v \ge 0, \quad \forall\, v\in \bar\calV \,\cap \,\calF_{z^*}\\
&\textstyle\sum_{a\in\bar\calS}\lambda_a a - \sum_{v\in\bar\calV\,\cap \,\calF_{z^*}}\pi_v v = \gamma d\\
\end{split}\tag{DLP}\label{eq:DLP}
\end{equation}
Since $z\in\conv(\bar\calS)$, $z^*\in \calF_{z^*} = \conv(\bar\calV\cap\calF_{z^*})$ (cf.\ Lemma~\ref{lem:conv_face}) and $\calZ$ is bounded, we know that:
\begin{enumerate}[label=\roman*)]
\item \eqref{eq:DLP} has an optimal solution, which we denote by $(\lambda^*,\pi^*,\gamma^*)$,
\item \eqref{eq:PLP} has an optimal solution, which we denote by $(p^*,t^*,\tau^*)$,
\item strong duality holds between~\eqref{eq:PLP} and~\eqref{eq:DLP}, i.e., $\gamma^*=t^*-\tau^*$.  
\end{enumerate} %and in addition, 

Let us make a few observations about~\eqref{eq:PLP} and~\eqref{eq:DLP}. 
First, %in~\eqref{eq:PLP}, 
we have $t^* = \max_{a\in\bar\calS}\, \lranglet{p^*}{a}$ and $
\tau^* = \min_{v\in\bar\calV\,\cap\,\calF_{z^*}}\, \lranglet{p^*}{v}$. In addition, since $z\ne z^*$, we have $\gamma^*>0$ and hence $t^*>\tau^*$. Also, since $\normtt{d}_{z^*}=1$, in~\eqref{eq:DLP}, we have 
\begin{align}
\textstyle \gamma^*= \normtt{u^*-w^*}_{z^*}, \quad \mbox{where}\quad u^*:= \sum_{a\in\bar\calS}\lambda^*_a a, \quad w^*:= \sum_{v\in\bar\calV\,\cap \,\calF_{z^*}}\pi^*_v v.  \label{eq:gamma*}
\end{align}
Next, from complementarity, we know that for all $a\in\bar\calS$, if $\lambda^*_a>0$, then $\ipt{p^*}{a} = t^*$, or equivalently, $a\in \bar\calS^*:= {\argmax}_{a\in\bar\calS}\, \lranglet{p^*}{a}\ne \emptyset$, and hence $u^*= \sum_{a\in\bar\calS^*}\lambda^*_a a\in\conv(\bar\calS^*)$.  %in other words,  %and therefore, $\{a\in\bar\calS:\lambda^*_a>0\}\subseteq \arg$ 
Similarly, we know that $w^*:= \sum_{v\in\bar\calV_{z^*}^*}\pi^*_v v \in \conv(\bar\calV_{z^*}^*)$, where $\bar\calV_{z^*}^*:=\argmin_{v\in\bar\calV\,\cap\,\calF_{z^*}}\, \lranglet{p^*}{v} \ne \emptyset$. 
From Lemma~\ref{lem:conv_argmin} and Lemma~\ref{lem:conv_face}, we know that 
\begin{equation}
\textstyle\emptyset\ne\conv(\bar\calV_{z^*}^*) = \argmin_{v\in\conv(\bar\calV\,\cap\,\calF_{z^*})}\, \lranglet{p^*}{v} = \argmin_{v\in\calF_{z^*}}\, \lranglet{p^*}{v}\in {\sf Faces} (\calF_{z^*}). \label{eq:conv_V*_face}
\end{equation}
Finally, since $t^*>\tau^*$ and $\tau^* = \ipt{p^*}{x}$ for all $x\in\conv(\bar\calV_{z^*}^*)$, we have $\bar\calS^*\cap \conv(\bar\calV_{z^*}^*)=\emptyset$. %, since %$\ipt{p^*}{a} = t^*$ for all $x\in\conv(\bar\calS^*)$ and  . %or equivalently, $\bar\calS^*\subseteq \bar\calV\setminus \bar\calV_{z^*}^*$. 
As such,  we have $ \bar\calS^*\subseteq \bar\calV\setminus \conv(\bar\calV_{z^*}^*)$ and hence 
\begin{equation}
\emptyset\ne \conv(\bar\calS^*)\subseteq \conv(\bar\calV\setminus \conv(\bar\calV_{z^*}^*))\quad\Longrightarrow\quad  \conv(\bar\calV_{z^*}^*)\ne \calZ. \label{eq:conv_barS*}
\end{equation}

%As a result, from~\eqref{eq:gamma*}, we have 
From the observations above, in particular~\eqref{eq:gamma*},~\eqref{eq:conv_barS*} and~\eqref{eq:conv_V*_face}, we have 
\begin{align*}
\Phi_\calS(x,\calX^*) =\gamma^* &= \normtt{u^*-w^*}_{z^*}\\
 &\ge {\min}_{u\in \conv(\bar\calS^*),\;w\in \conv(\bar\calV_{z^*}^*)} \;\normtt{u-w}_{z^*}\\
&\ge {\min}_{u\in \conv(\bar\calV\setminus \conv(\bar\calV_{z^*}^*)),\;w\in \conv(\bar\calV_{z^*}^*)}\; \normtt{u-w}_{z^*}\\
&\ge {\min}_{\emptyset\ne \calF\in {\sf Faces} (\calF_{z^*}), \,\calF\ne \calZ}\;{\min}_{u\in \conv(\bar\calV\setminus \calF),\;w\in \calF} \;\normtt{u-w}_{z^*}\\
&= {\min}_{\emptyset\ne \calF\in {\sf Faces} (\calF_{z^*}), \,\calF\ne \calZ}\; \dist_{\normtt{\cdot}_{z^*}}(\conv(\bar\calV\setminus \calF), \calF). \nt\label{eq:min_F_dist}
 %:= \dist_{\normtt{\cdot}_{z^*}} \big(\conv(\bar\calS^*),\conv(\bar\calV_{z^*}^*)\big). \
\end{align*}
Note that the RHS of~\eqref{eq:min_F_dist} only depends on $\bar\calV$ and $z^*$ (and does not depend on $x$). 
Therefore, by taking minimum on $\Phi_\calS(x,\calX^*)$ over $\calS\in{\sf AR}(x)$ and $x\in\calX\setminus\calX^*$, we complete the proof. \qed \vspace{2ex}

%As such, by taking minimum on the LHS of~\eqref{eq:min_F_dist} over $\calS\in\$
%for all $v\in\bar\calV\cap\calF_{z^*}$, if $\pi^*_v>0$, then $v\in \bar\calV^*:=\argmin_{v\in\bar\calV\,\cap\,\calF_{z^*}}\, \lranglet{p^*}{v}$, and . 

%\newpage

\noindent {\em Proof of Lemma~\ref{lem:conv_face}.}  
Note that since $\calV\cup\calF\subseteq \calF$ and $\calF$ is convex, we clearly have $\conv(\calV\cup\calF)\subseteq \calF$. To show $\calF\subseteq \conv(\calV\cup\calF)$, it suffices to only consider the case that $\calF\ne \emptyset$ (otherwise the statement trivially holds).  Fix any $x\in\calF$ and consider its ``active atoms'' in $\calV$, namely $\calV'\subseteq\calV$ such that $x\in\ri \conv(\calV)$. Since $\calF$ is a face of $\conv(\calV)$, we have $\calV'\subseteq \calF$ and hence $\calV'\subseteq \calV\cup\calF$. Then we know that $x\in\conv(\calV')\subseteq \conv(\calV\cup\calF)$, and hence $\calF\subseteq \conv(\calV\cup\calF)$. \qed \vspace{2ex}

\noindent {\em Proof of Lemma~\ref{lem:conv_argmin}.}  
For convenience, let $\calV^*:= \argmin_{v\in\calV}\, \ipt{p}{v}$, $\calX:= \conv(\calV)$ and $\calX^*:= \argmin_{x\in\calX}\, \ipt{p}{x}$. We show that for any $x\in\calX^*$ and any $v\in\calV^*$, $\ipt{p}{x} = \ipt{p}{v}$. Indeed, let $x = \sum_{v'\in\calV}\, \lambda_{v'} v'$, where $\sum_{v'\in\calV}\, \lambda_{v'}=1$ and $\lambda_{v'}\ge 0$ for all $v'\in\calV$, and we have
\begin{equation}
\ipt{p}{v}\ge \ipt{p}{x} = \textstyle\sum_{v'\in\calV}\, \lambda_{v'} \ipt{p}{v'}\ge \ipt{p}{v}. 
\end{equation}
In addition, if $\lambda_{v'} >0$, then $\ipt{p}{v'}=\ipt{p}{v}$, and hence $v'\in\calV^*$. Consequently, $x\in\conv(\calV^*)$ and hence $\calX^*\subseteq \conv(\calV^*)$. On the other hand, since $\ipt{p}{v} = \ipt{p}{x}$, we have $v\in\calX^*$ and hence $\calV^*\subseteq \calX^*$. Since $\calX^*$ is convex, we have $\conv(\calV^*)\subseteq \calX^*$. 
\qed

%\section{Details of Remark~\ref{rmk:minimal_face} }
\section{Proof of the Minimality of $\calF^*$}\label{app:details_rmk}

%We shall show that $\calF^*$
%In the following, 
We use the same notations in the proof of Lemma~\ref{lem:exist_face}. In addition, for any set $\calS\neq \emptyset$, denote %its (affine) dimension as $\dim \calS$ and 
the linear space associated with $\aff\calS$ as $\lin\calS$ (namely, $\lin\calS:= \aff\calS-x$ for any $x\in\aff\calS$).  Also, denote the nullspace of $\rvE$ as ${\sf null}(\rvE)$.

Let  $-\nabla \barF(z^*)\in\ri \calN_{\calZ}(z^*) = \{\sum_{i\in\calI^*}\,\gamma_i c_i: \gamma_i>0,\;\forall\,i\in\calI^*\}$, and hence $\ipt{\nabla \barF(z^*)}{z - z^*} = 0$ %\, \Longleftrightarrow \, 
if and only if $z\in\calF_{\calI^*}$. Note that $\calF_{\calI^*}$ is the minimal face of $\calZ$ that contains $z^*$. 
Let $\calF^*:= \{x\in\calX:\rvE x \in \calF_{\calI^*}\}$, and  we claim that $\calF^*$ is the minimal face of $\calX$ that contains $\calX^*$. %Suppose the opposite is true. 
%Suppose this is not true. 
If not, let $\calG$ be the minimal proper face of $\calF^*$ that contains $\calX^*$, so that %i) there exists a vertex of $\calF^*$, say $v$, such that $v\not\in\aff\calG$, and ii) 
there exists $x^*\in\calX^*$ such that $x^*\in\ri\calG$. Denote the set of vertices of %$\calF^*$ and 
of $\calG$ as %$\{v_j\}_{j=1}^p$ and 
$\{v_j\}_{j\in\calJ}$. %, respectively, where $\calJ\subsetneqq [p]$. 
%$\calG$ by $\{v_j\}_{j=1}^p$, and 
Since $x^*\in\ri\calG$, we know that $x^* = \sum_{j\in\calJ} \alpha_j v_j $ for some $\{\alpha_j\}_{j\in\calJ}\subseteq\bbR_{++}$ and $\sum_{j\in\calJ} \alpha_j = 1$. % and 
  %then exists a proper face of $\calF^*$, say $\calG$, that contains $\calX^*$. 

%Since $\calF^*$ is a polytope, there exist a vertex of $\calF^*$, say $v$, such that $v\not\in\aff\calG$. Now,  note that 
Fix any $v\in \calF^*\setminus \calG$.  
 %We first show that 
 Let us note that ${\sf null}(\rvE)\cap\lin(\calG\cup\{v\})={\sf null}(\rvE)\cap\lin(\{v_j\}_{j\in\calJ}\cup\{v\}) \subseteq \lin\calG$. 
 If not, there exists $d \in{\sf null}(\rvE)$ and $d=\beta v + \sum_{j=1}^p \beta_j v_j$ 
 for some $\beta\ge 0$ and $\beta + \sum_{j=1}^p \beta_j=0$ such that $d\not\in\lin\calG$. %As a result, 
Since $\{\alpha_j\}_{j\in\calJ}\subseteq\bbR_{++}$, for sufficiently small $\varepsilon>0$, we have 
\begin{equation}
x^*+\varepsilon d = \textstyle \sum_{j\in\calJ} (\alpha_j + \varepsilon \beta_j)v_j + \varepsilon\beta v \in \conv(\{v_j\}_{j\in\calJ}\cup\{v\})\subseteq \calF^*\subseteq \calX.  
\end{equation}
By the definition of $\calX^*$ (cf.~\eqref{eq:calX*}), we know that $x^*+\varepsilon d\in\calX^*\subseteq\aff\calG$, and hence $d\in\lin\calG$, leading to a contradiction. % However, since $d\not\in\lin\calG$, we have $x^*+\varepsilon d\not\in\$
%we know that Clearly, 
% such that $d \in\calN(\rvE)$. 
%Next, we show that 
Based on this, we have $\rvE v\not\in \rvE(\calG)$. %Otherwise,
If not, we have $\rvE v = \sum_{j\in\calJ} \baralpha_j \rvE v_j$ for some $\{\baralpha_j\}_{j\in\calJ}\subseteq\bbR_{+}$ and $\sum_{j\in\calJ} \baralpha_j = 1$, %which implies that 
and hence $\barv:=v-\sum_{j\in\calJ} \baralpha_j v_j \in{\sf null}(\rvE)$. Since $\barv\in\lin(\calG\cup\{v\})$, we have $\barv\in\lin\calG$, and hence $v\in\aff\calG$. However, since $v\in \calF^*\setminus \calG$, $v\not\in\aff\calG$ -- this is a contradiction. 

Next, we show that $\rvE(\calG)$ is a proper face of $\rvE(\calF^*) = \calF_{\calI^*}$. For any $u,w\in\calF_{\calI^*}$, there exists $x,y\in\calF^*$ such that $u = \rvE x$ and $w = \rvE y$. Suppose $\lambda u+(1-\lambda)w=\rvE(\lambda x+(1-\lambda)y)\in\rvE(\calG)$, since $\lambda x+(1-\lambda)y\in\calF^*$, we must have $\lambda x+(1-\lambda)y\in\calG$ (otherwise, let $v=\lambda x+(1-\lambda)y$ in the above and we have a contradiction). Since $\calG$ is a face of $\calF^*$, we have $x,y\in\calG$ and hence $u,w\in\rvE(\calG)$. In addition, since there exists $v\in\calF^*\setminus \calG$ such that $\rvE v\not\in\rvE(\calG)$, $\rvE(\calG)\subsetneqq \calF_{\calI^*}$. 

Now, since $z^* = \rvE x^*\in\rvE(\calG)$, we know that   $\calF_{\calI^*}$ is not the minimal face of $\calZ$ that contains $z^*$, leading to a contradiction. 
%If this is true

%Indeed, if  $\rvE v \rvE(\calG)$

\section{Proof of Lemma~\ref{lem:delta_k_ub} } \label{app:proof_delta_k_ub}

The  proof of Lemma~\ref{lem:delta_k_ub} leverages the following lemma.

\begin{lemma} \label{lem:abs_ip_diff}
In Algorithm~\ref{algo:AFW_SC}, if $\normt{y^k-y^*}_{y^*}< 1$, then %for all $x^*\in\calX^*$,  %for any $x\in\calX$, we have %$k\ge 0$,
\begin{equation}
{\sup}_{x^*\in\calX^*}\;{\sup}_{x\in\calX}\;\abst{\ipt{\nabla F(x^k)}{x-x^k} - \ipt{\nabla F(x^*)}{x-x^*}}\le \frac{\normt{y^k-y^*}_{y^*} }{1-\normt{y^k-y^*}_{y^*}}R({y^*}) + G_k, %\quad  \forall\,x\in\calX, 
\label{eq:abs_ip_diff}
\end{equation} 
where $R(y^*)$ is defined in~\eqref{eq:def_R_y*}. 
\end{lemma}

\begin{proof}
Fix any $x\in\calX$ and any $x^*\in\calX^*$. Let $y:=\rvA x\in\calY$,  % and recall that $y^*= \rvA x^*$, for all $x^*$. 
and we have 
\begin{align}
&\abst{\ipt{\nabla F(x^k)}{x-x^k} - \ipt{\nabla F(x^*)}{x-x^*}}\\
=\;&\abst{\ipt{\nabla F(x^k)-\nabla F(x^*)}{x-x^*} + \ipt{\nabla F(x^k)}{x^*-x^k}}\\
\le \; &\abst{\ipt{\nabla f(y^k)-\nabla f(y^*)}{y-y^*}} + \ipt{\nabla F(x^k)}{x^k-x^*}\\
\le \; & \normt{\nabla f(y^k)-\nabla f(y^*)}_{y^*}^*\normt{{y-y^*}}_{y^*} + G_k,
\end{align}
where the last step follows from% the definition of $G_k$ in
~\eqref{eq:def_G_k}. Now, using Lemma~\ref{lem:local_norm_grad_diff}, we complete the proof. 
\end{proof}

Now, note that from~\eqref{eq:barF_lb2}, we see that $\normt{y^k-y^*}_{y^*}\le \sqrt{\delta_k}/\sqrt{\mu_{R(y^*)}}$, and from~\eqref{eq:delta_k_le}, we have %$$. Hence
\begin{equation}
\normt{y^k-y^*}_{y^*}< \frac{\Delta}{2R(y^*)+\Delta}<1 \;\;\;\Longrightarrow\;\;\; \frac{\normt{y^k-y^*}_{y^*} }{1-\normt{y^k-y^*}_{y^*}}R({y^*})< \frac{\Delta}{2R(y^*)}R(y^*) = \frac{\Delta}{2}.  \label{eq:y^k-y^*}
\end{equation}
In addition, from the definition of $D_k$ in Step~\ref{item:step_size_SC_a}\ref{item:adapt_step} and that $\normt{y-y^*}_{y^*}<1$, we have 
\begin{align}
D_k\le {\sup}_{y\in\calY}\; \normt{y-y^k}_{y^k} \lea {\sup}_{y\in\calY}\;\frac{\normt{y-y^*}_{y^*}+\normt{y^k-y^*}_{y^*}}{1-\normt{y^k-y^*}_{y^*}} = \frac{R(y^*)+\normt{y^k-y^*}_{y^*}}{1-\normt{y^k-y^*}_{y^*}}, \label{eq:D_k_ub2} %\le \frac{(2R(y^*)+\Delta)R(y^*) + \Delta}{2R(y^*)}. 
%&\le {\sup}_{y\in\calY}\; \normt{y-y^*}_{y^k} +\normt{y^k-y^*}_{y^k} \\
%&\le {\sup}_{y\in\calY}\; \frac{\normt{y-y^*}_{y^*}}{1-\normt{y^k-y^*}_{y^k}} +\normt{y^k-y^*}_{y^k}
\end{align}
where (a) follows from Lemma~\ref{lem:approx_Hessian}. 
%Therefore, 
From~\eqref{eq:r_k_ub2},%
~\eqref{eq:D_k_ub2},~\eqref{eq:y^k-y^*} and~\eqref{eq:delta_k_le}, we have $\delta_k<1/16$ and 
\begin{align*}
G_k\le 4\max\{\delta_k,\sqrt{\delta_k}\}D_k = 4\sqrt{\delta_k}D_k< 4\sqrt{\delta_k}\frac{R(y^*)+{\Delta}/({2R(y^*)+\Delta})}{2R(y^*)/({2R(y^*)+\Delta})} < \frac{\Delta}{2}. 
\end{align*}
%hence $$. % (since $\delta_k<1/16$). 
%\begin{equation}
%G_k\le 4\max\{\delta_k,\sqrt{\delta_k}\}D_k\le 4\sqrt{\delta_k}\frac{(2R(y^*)+\Delta)R(y^*) + \Delta}{2R(y^*)}\le \frac{\Delta}{2}. 
%\end{equation}
%from, we have
Now, using~\eqref{eq:abs_ip_diff}, % in Lemma~\ref{lem:abs_ip_diff}, 
we complete the proof.

\bibliographystyle{spmpsci}      % mathematics and physical sciences

\bibliography{GF-papers-nips-better,math_opt,stat_ref,mach_learn}

\end{document}